\documentclass[12pt]{article}
\usepackage[mathscr]{eucal}
\usepackage{amssymb}
\usepackage{latexsym}
\usepackage{theorem}
\theoremstyle{plain}
\newtheorem{Theorem}{Theorem}[section]

\newtheorem{Lemma}[Theorem]{Lemma}
\newtheorem{Proposition}[Theorem]{Proposition}

\theorembodyfont{\rmfamily}
\newtheorem{Remark}[Theorem]{Remark}
\newtheorem{Example}[Theorem]{Example}
\newtheorem{Definition}[Theorem]{Definition}
\newtheorem{pf}{Proof}

\begin{document}

\begin{titlepage}
Differential Geometry and Integrable Systems,\\
{\it Contemporary Mathematics AMS}, vol.308, 2002

\vskip 1cm

\begin{center}
{\bf {\Large Weierstra{\ss} type representation of timelike
surfaces with constant mean curvature}}
\end{center}

\vskip 2cm

\begin{center}
\bf{Josef Dorfmeister, Junichi Inoguchi, Magdalena Toda}
\end{center}

\vskip 2cm

{\abstract We derive a correspondence between (Lorentzian)
harmonic maps into the pseudosphere $S_1^2$, with appropriate
regularity conditions, and certain connection 1-forms. To these
harmonic maps, we associate a representation of type Weierstrass,
and we apply it to construct timelike surfaces with constant mean
curvature.}

\date{}

\end{titlepage}

\section*{Introduction}

As is well known, every solution to the elliptic sinh-Gordon (or
sinh-Laplace) equation:
$$
\omega_{z {\bar z}}+\sinh \omega=0
$$
describes a constant mean curvature (CMC) surface in Euclidean
$3$-space $\mathbb{E}^3$. The symmetric quadratic form
$I=e^{\omega}dzd{\bar z}$ is the induced metric of the CMC
surface. With respect to the conformal structure determined by
$I$, the Gau{\ss} map is a harmonic map into $S^2$.

F.~Pedit, H.~Wu and the first named author established a loop
group theoretic Weierstra{\ss}-type representation for harmonic
maps of Riemann surfaces into compact Riemannian symmetric spaces
\cite{DPW}. This loop group theoretic Weierstra{\ss}-type
representation is frequently referred to in the literature as the
{\it DPW-method}.

Furthermore G.~Haak and the first named author used the DPW-method
intensively for the investigation of
 CMC surfaces \cite{DH}.
H.~Wu gave a simple way for determining a normal form of the
Weierestra{\ss} data for CMC surfaces in \cite{Wu}. In \cite{DPT},
Pedit, the first and third named author gave a reinterpretation of
the classical Weierstra{\ss} representation of minimal surfaces in
terms of the DPW-method.

\vspace{0.2cm}

On the other hand, every solution $\phi$ to the sine-Gordon
equation:
$$
\phi_{xy}+\sin \phi=0.
$$
describes a (weakly regular) pseudospherical surface in Euclidean
$3$-space $\mathbb{E}^3$. The solution $\phi$ is the angle of two
asymptotic directions. The quadratic form $II=\sin \phi \>{\mathrm
d}x {\mathrm d}y$ is the second fundamental form of the surface.
With respect to the Lorentzian conformal structure determined by
$II$, the Gau{\ss} map is a (Lorentzian) harmonic map into the
unit $2$-sphere $S^2$.

M.~Melko and I.~Sterling \cite{MS}-\cite{MS2} presented a modern
approach to pseudospherical surfaces via the theory of finite type
(Lorentzian) harmonic maps into $S^2$.

Recently the third named author established a Weierstra{\ss}-type
representation for pseudospherical surfaces in $\mathbb{E}^3$ in
her thesis \cite{T}. In particular, she showed that there exist
certain normalized potentials for pseudospherical surfaces.
Finally, in \cite{DS} it was shown that finite type
pseudospherical surfaces can be constructed from particularly
simple normalized potentials.

\vspace{0.2cm}

Both, the  sine-Gordon and the elliptic sinh-Gordon equation are
real forms of the complexified sine-Gordon equation.
 But the complexified sine-Gordon equation
has still different  real forms: the hyperbolic sinh-Gordon and
the hyperbolic cosh-Gordon equations. Solutions to these equations
do not describe CMC surfaces in $\mathbb{E}^3$ but in Minkowski
$3$-space $\mathbb{E}^3_1$.

\vspace{0.2cm} In fact, let $M$ be a timelike CMC surface in
$\mathbb{E}^3_1$ parametrized locally by null coordinates $(x,y)$.
Denote  by $D$ the discriminant of the characteristic equation for
the shape operator. Then the Gau{\ss}-Codazzi equations of $M$
become
$$
\omega_{xy}+H^2\sinh \omega=0, \ \mathrm{if} \ D>0,
$$
$$
\omega_{xy}+\frac{H^2}{2}e^{\omega}=0, \ \mathrm{if} \ D=0,
$$
$$
\omega_{xy}+H^2\cosh \omega=0, \ \mathrm{if} \ D<0.
$$

\smallskip

At this point  we  would like to point out a certain similarity
between timelike CMC surfaces and CMC surfaces in hyperbolic
$3$-space. Note that both $H^3$ and $\mathbb{E}^3_1$ are naturally
imbedded in Minkowski $4$-space $\mathbb{E}^4_1$.

Let $M$ be a CMC surface in hyperbolic $3$-space $H^3$
parametrized by isothermic coordinates (isothermal-curvature line
coordinates) $z=x+\sqrt{-1}y$, then its Gau{\ss}-Codazzi equations
become
$$
\omega_{z{\bar z}}+(H^2-1)\sinh \omega=0, \ \mathrm{if} \ H^2>1,
$$
$$
\omega_{z{\bar z}}-\frac{1}{2}e^{-\omega}=0, \ \mathrm{if} \
H^2=1,
$$
$$
\omega_{z{\bar z}}+(H^2-1)\cosh \omega=0, \ \mathrm{if} \ H^2<1.
$$

Therefore, at the level of Gau{\ss}-Codazzi equations, timelike
CMC surface geometry can be considered as a ``hyperbolic version"
of CMC surface geometry in $H^3$. There is another similarity
between timelike surface geometry in $\mathbb{E}^3_1$ and surface
geometry in $H^3$. In fact, timelike HIMC surfaces, {\it i.e.},
timelike surfaces with harmonic inverse mean curvature in
$\mathbb{E}^3_1$ \cite{FI} and Bonnet surfaces in $H^3$, (which
are not Willmore surfaces) \cite{BE} are described by the same
integrable equation, namely the  {\it Painlev{\'e} equations of
type} V and VI.

\vspace{0.2cm}

The hyperbolic sinh-Gordon equation and the  Liouville equation
have been studied extensively by the soliton theoretic approach.
For instance, Babelon and Bernard \cite{BB} studied (the
infinitesimal version of) the dressing transformations for the
hyperbolic sinh-Gordon equation. L.~McNertney studied B{\"a}cklund
transformations for the hyperbolic sinh-Gordon and the Liouville
equation via the classical theory of line-congruences in her
thesis \cite{Mc}. H.-S.~Hu \cite{Hu2} and the second named author
\cite{I2} gave a Darboux form of B{\"a}cklund transformations on
the hyperbolic sinh-Gordon equation. Namely, in \cite{Hu2} and
\cite{I2}, B{\"a}cklund transformations are reformulated as {\it
transformations on extended framings}.

\vspace{0.2cm}

 However as far as the authors know,
only few studies on the cosh-Gordon equation are available.
M.~V.~Babich obtained finite-gap solutions to the elliptic
cosh-Gordon equation \cite{Ba}. Babich and A.~I.~Bobenko studied
minimal surfaces in $H^3$ in terms of finite-gap solutions of the
elliptic cosh-Gordon equation \cite{BaBo}. V.~Y.~Novokshenov
studied radial-symmetric solutions to the elliptic cosh-Gordon
equation \cite{Nov2}. The radial-symmetry reduces the elliptic
cosh-Gordon equation to the third Painlev{\'e} equation. Moreover
he studied the minimal surfaces in $H^3$ corresponding to these
solutions to the elliptic cosh-Gordon equation \cite{Nov}.

In this paper we establish a Weierstra{\ss}-type representation
for timelike CMC surfaces in Minkowski $3$-space. The
Weierstra{\ss}-type representation gives a unified theory of
constructing solutions to the sinh-Gordon, the Liouville and the
cosh-Gordon equations. Moreover our Weierstra{\ss}-type
representation is regarded as {\it nonlinear d'Alembert formula}
for these three nonlinear wave equations.

\vspace{0.2cm}

This paper is organized as follows:

After establishing the requisite facts on geometry of surfaces in
Minkowsi 3-space in Section 1, we devote Section 2 to prepare
ingredients from loop group theory.

In Section 3, we derive a correspondence between harmonic maps
into the pseudosphere $S^2_1$ with appropriate regularity and a
certain kind of connection one-forms.

The Weierstra{\ss}-type representation for (Lorentzian) harmonic
maps into $S^2_1$ is introduced in Section 4. We apply the
Weierstra{\ss}-type representation for constructing timelike CMC
surfaces in Section 5.

In the final section, we discuss fundamental examples of timelike
CMC surfaces via the Weierstra{\ss}-type representation.

\section{Timelike surfaces}
\noindent {\large {\bf 1.1}} \hspace{0.15mm} We start with
preliminaries on the geometry of timelike surfaces in Minkowski
3-space. \par

Let $\mathbb{E}^3_1$ be {\it Minkowski} $3$-{\it space} with
Lorentzian metric $\langle \cdot ,\cdot  \rangle$. The metric
$\langle \cdot  ,\cdot  \rangle$ is expressed as $\langle \cdot
,\cdot  \rangle= -{\mathrm d}u^2_1+ {\mathrm d}u^2_2+{\mathrm
d}u^2_3$ in terms of the natural coordinate system $(u_1,u_2,u_3)$
of the Cartesian $3$-space $\mathbb{R}^3$.

Let $M$ be a connected orientable $2$-manifold and $\varphi : M
\rightarrow \mathbb{E}^3_1$ an immersion. The immersion $\varphi$
is said to be {\it timelike} if the induced metric $I$ of $M$ is
Lorentzian. The induced Lorentzian metric $I$ determines a Lorentz
conformal structure ${\mathcal C}_I$ on $M$. We treat
$(M,{\mathcal C}_I)$ as a Lorentz surface and $\varphi$ as a
conformal immersion. For the general theory of Lorentz surfaces,
we refer to T.~Weinstein \cite{We}.

\vspace{0.2cm}

Hereafter we will assume that $M$ is an orientable timelike
surface in ${\mathbb E}^3_1$ (immersed by $\varphi$).

It is worthwhile to remark that there exists no compact timelike
surface in $\mathbb{E}^3_1$. (See B.~O'Neill \cite{ON}, p. 125.)

\vspace{0.2cm}

Let $(x,y)$ be a {\it null coordinate system} with respect to the
conformal structure ${\mathcal C}_{I}$. Then the first fundamental
form $I$ is written in terms of $(x,y)$ as follows:
$$
I=e^{\omega}\>{\mathrm d}x{\mathrm d}y.
$$

Now let $N$ be a unit normal vector field of $M$. Namely a vector
field $N$ along $M$ satisfying
$$
\langle N, N \rangle=1,\ \langle \varphi_x,N \rangle= \langle
\varphi_y,N \rangle=0.
$$

The {\it second fundamental form} $II$ of $M$ derived from $N$ is
defined by
$$
II=-\langle {\mathrm d}\varphi,\  {\mathrm d}N \rangle.
$$
The {\it shape operator} $S$ of $M$ derived from $N$ is
$$
S:=-{\mathrm d}N.
$$
The shape operator $S$ is related to $II$ by
$$
II(X,Y)=\langle SX,Y \rangle
$$
for all vector fields $X,\ Y$ on $M$. The mean curvature $H$ of
$M$ is defined by
$$
H=\frac{1}{2}\mathrm{tr}\>S.
$$
Note that $H$ is computed by the following formula:
$$
H=\frac{1}{2}\mathrm{tr}(II \cdot I^{-1}).
$$
Note that the Gau{\ss}ian curvature $K$ of $M$ is computed as
$$
K=\det S=\det (II \cdot I^{-1}).
$$
(See \cite{ON}, p.~107.) The characteristic values of $S$, {\it
i.e.,} the (complex) solutions to
$$
\det (t\>{\mathrm I}-S)=0,\ \mathrm{I}=\mathrm{identity \ of \
}TM.
$$
are called the {\it principal curvatures}. Since the metric $I$ is
indefinite, both principal curvatures may be {\it non real}
complex numbers. It is easy to check from the definitions that $H$
is the mean of the two principal curvatures and $K$ is the product
of the two principal curvatures.

A point $p$ of $M$ is said to be an {\it umbilic point} if $II$ is
proportional to $I$ at $p$. Equivalently, $p$ is an umbilic point
if and only if the two principal curvatures at $p$ are the same
real number and the corresponding eigenspace is $2$-dimensional.

A timelike surface is said to be a {\it totally umbilic surface}
if all the points are umbilical.

It is known that every totally umbilic timelike surface in
$\mathbb{E}^3_1$ is congruent to an open portion of a pseudosphere
$$
S^2_1(r):= \left \{ {\mathbf u} \in \mathbb{E}^3_1\ \vert \
\langle \mathbf{u},\mathbf{u} \rangle =r^2 \right \}
$$
of radius $r>0$ or a timelike plane.
\medskip

Let $\varphi:M \to \mathbb{E}^3_1$ be a timelike surface with unit
normal vector field $N$ as before. Then, on a simply connected
null coordinate region ${\mathbb D}$, we can define an orthonormal
frame field $\mathcal{F}$ defined by
$$
{\mathcal F}=(e^{-\omega/2}(-\varphi_{x}+\varphi_{y}),
e^{-\omega/2} (\varphi_{x}+\varphi_{y}), N):{\mathbb D}\to
\mathrm{O}^{++}_{1}(3),
$$
where $\mathrm{O}^{++}_{1}(3)$ denotes the identity component of
the Lorentz group
\[
\mathrm{O}_{1}(3)= \left \{ A \in \mathrm{GL}(3;{\mathbb R})\
\vert \ \langle A \mathbf{u},A\mathbf{v} \rangle=\langle
\mathbf{u},\mathbf{v} \rangle,\ \mathbf{u}, \mathbf{v} \in
\mathbb{E}^3_1 \ \right \}.
\]

\vspace{0.2cm}

\noindent {\large {\bf 1.2}}\hspace{0.15cm} Throughout this paper,
we identify ${\mathbb E}^3_1$ with the Lie algebra ${\mathfrak
g}=\mathfrak{sl}(2;{\mathbb R})$. We take the following basis
$\{\>{\mathbf i}, \> {\mathbf j}^{\prime},\> {\mathbf
k}^{\prime}\}$ of $\mathfrak{sl}(2;{\mathbb R})$:

\[
{\mathbf i}= \left (
\begin{array}{cc}
0 & -1 \\
1 & 0
\end{array}
\right ), \ \ {\mathbf j}^{\prime} =\left (
\begin{array}{cc}
0 & 1 \\
1 & 0
\end{array}
\right ),\ \ {\mathbf k}^{\prime} =\left (
\begin{array}{cc}
-1 & 0 \\
0 & 1
\end{array}
\right ).
\]
The basis $\{\>{\mathbf i}, \> {\mathbf j}^{\prime},\> {\mathbf
k}^{\prime}\}$ satisfies the following relation:
$$
{\mathbf i}^2=-{\mathbf 1},\ {\mathbf j}^{\prime 2}= {\mathbf
k}^{\prime 2}={\mathbf 1},
$$

$$
{\mathbf i}{\mathbf j}^{\prime}= -{\mathbf j}^{\prime} {\mathbf
i}= {\mathbf k}^{\prime},\ {\mathbf j}^{\prime} {\mathbf
k}^{\prime}= -{\mathbf k}^{\prime} {\mathbf j}^{\prime}= -{\mathbf
i},\ {\mathbf k}^{\prime}{\mathbf i} =-{\mathbf i}{\mathbf
k}^{\prime} ={\mathbf j}^{\prime}.
$$
Here ${\mathbf 1}$ denotes the identity matrix:
\[
\mathbf{1}= \left (
\begin{array}{cc}
1 & 0 \\
0 & 1
\end{array}
\right ).
\]

Hereafter we identify $\mathbb{E}^3_1$ with $\mathfrak{g}$ via
this basis.
$$
(u_1,u_2,u_3) \longleftrightarrow u_1{\mathbf i}+ u_2{\mathbf
j}^{\prime} +u_3{\mathbf k}^{\prime}. \eqno(1.1)
$$

The real algebra $\mathbb{H}^{\prime}$ generated by $\{\> {\mathbf
1},\> {\mathbf i}, \> {\mathbf j}^{\prime},\> {\mathbf
k}^{\prime}\}$ is called the algebra of {\it split quaternions}.
The algebra $\mathbb{H}^\prime$ is isomorphic to the algebra
$\mathrm{M}(2;\mathbb{R})$ of all 2 by 2 real matrices. The
commutation relations of $\mathfrak{g}$ are given by

$$
[\ {\mathbf i},\> {\mathbf j}'\ ]=2\> {\mathbf k}',\ [\ {\mathbf
j}',\> {\mathbf k}'\ ]=-2\> {\mathbf i},\ [\ {\mathbf k}',\>
{\mathbf i}\ ]=2\> {\mathbf j}'.
$$

\vspace{0.2cm}

By the linear isomorphism (1.1) the Lorentz metric $\langle \cdot,
\cdot \rangle$ corresponds to the scalar product:
$$
\langle X,Y \rangle= \frac{1}{2}\mathrm{tr}(XY),\ \ X,Y \in
\mathfrak{g}. \eqno(1.2)
$$

This scalar product induces a biinvariant Lorentz metric of
constant curvature $-1$ on the special linear group
$G=\mathrm{SL}(2;\mathbb{R})$. Hence $G$ is identified with the
anti-de Sitter $3$-space $H^3_1$. (See \cite{ON}.)

\vspace{0.2cm}

The special linear group $G$ acts isometrically on $\mathfrak{g}$
via the $\mathrm{Ad}$-action:
$$
\mathrm{Ad}:G \times \mathfrak{g}\to \mathfrak{g};\ \
\mathrm{Ad}(g)X=gXg^{-1}.
$$
The Ad-action induces a double covering $G \to
\mathrm{O}^{++}_{1}(3)$ of the Lorentz group
$\mathrm{O}^{++}_{1}(3)$.

\vspace{0.2cm}

By using this double covering we can find a  lift ${\hat \Phi}$
(called a {\it coordinate frame}) of $\mathcal{F}$ to
$\mathrm{SL}(2;{\mathbb R})$:
$$
\mathrm{Ad}({\hat \Phi}) ({\mathbf i},{\mathbf j}^{\prime},
{\mathbf k}^{\prime})= \mathcal{F}. \eqno(1.3)
$$
The coordinate frame ${\hat \Phi}$ satisfies the following Frenet
( or Gauss-Weingarten) equations \cite{I}:

$$
\frac{\partial}{\partial x}{\hat \Phi} ={\hat \Phi}{\hat U},\
\frac{\partial}{\partial y} {\hat \Phi}={\hat \Phi}{\hat V},
\eqno(1.4)
$$
where

$$
{\hat U}= \left (
\begin{array}{cc}
-\frac{1}{4}\omega_{x}
& -Qe^{-\frac{\omega}{2}} \\
\frac{H}{2}e^{\frac{\omega}{2}} & \frac{1}{4}\omega_x
\end{array}
\right ), \ \ {\hat V}=\left (
\begin{array}{cc}
\frac{1}{4}\omega_y & -\frac{H}{2}e^{\frac{\omega}{2}} \\
R e^{-\frac{\omega}{2}} & -\frac{1}{4}\omega_y
\end{array}
\right ), \eqno(1.5)
$$
and $Q:=\langle \varphi_{xx},N \rangle$, $R:=\langle
\varphi_{yy},N \rangle$, $H=2e^{-\omega} \langle \varphi_ {xy},N
\rangle $.

The function $H$ coincides with the mean curvature of $\varphi$.
It is easy to see that $Q^{\#}:=Q{\mathrm d}x^2$ and
$R^{\#}:=R{\mathrm d}y^2$ are globally defined on $M$. The
quadratic differentials  $Q^{\#}$ and $R^{\#}$ are called the {\it
Hopf differentials of \/} $M$. \par The second fundamental form
$II$ of $M$ is related to $Q$ and $R$ by
$$
II=Q^{\#}+R^{\#}+HI.
$$
This formula implies that the common zeros of $Q$ and $R$ coincide
with the umbilic points of $M$.

The Gauss equation which describes a relation between $K$, $H$ and
$Q$ takes the following form:
$$
H^2-K=4e^{-2 \omega}Q R.
$$

Note that the condition $QR=0$ does not imply the umbilicity of
$M$. (See T.~K.~Milnor \cite{Mil}).
\bigskip

The Gauss-Codazzi equation, {\it i.e.}, the integrability
condition of the Frenet equations,
$$
{\hat V}_{x}- {\hat U}_{y}+[\ {\hat U},{\hat V}\ ]=0
$$
has the following form:
$$
\omega_{xy}+\frac{1}{2}H^{2}e^{\omega}-2QRe^{-\omega}=0,
\eqno{(\mathrm{G})}
$$

$$
H_{x}=2e^{-\omega}Q_y,\ \ H_{y}=2e^{-\omega}R_x.
\eqno{(\mathrm{C})}
$$

The Codazzi equations (C) show that the constancy of the mean
curvature $H$ is equivalent to the condition $Q_y=R_x=0$, {\it
i.e.,\/} $Q=Q(x),\ R=R(y)$.

\begin{Remark}
Let $(M,\mathcal{C})$ be a Lorentz surface and $(x,y)$ a null
coordinate system. Then the following differential operators are
well defined:
$$
{\mathrm d}^{\prime}:=\frac{\partial}{\partial x}{\mathrm d}x,\ \
{\mathrm d}^{\prime \prime}:= \frac{\partial}{\partial y}{\mathrm
d}y. \eqno(1.6)
$$
A function $f:M \to \mathbb{R}$ is said to be a {\it Lorentz
holomorphic function} [resp. {\it Lorentz anti-holomorphic
function}] if ${\mathrm d}^{\prime \prime}f=0$ [resp. ${\mathrm
d}^{\prime}f=0$].

Next a $1$-form $A=A_x {\mathrm d}x+A_y {\mathrm d}y$ is said to
be a {\it Lorentz holomorphic} $1$-{\it form} if $A=A_x {\mathrm
d}x$ and $A_x$ is a Lorentz holomorphic function. Similarly $A$ is
said to be a {\it Lorentz anti-holomorphic} $1$-{\it form} if
$A=A_y {\mathrm d}y$ and $A_y$ is a Lorentz anti-holomorphic
function.

According to these terminologies, the constancy of mean curvature
is characterized as follows:

\vspace{0.2cm}

{\it Let} $\varphi:M \to \mathbb{E}^3_1$ {\it be a timelike
surface. Then $(M,\varphi)$ is of constant mean curvature if and
only if $Q$ is a Lorentz holomorphic function and $R$ is a Lorentz
anti-holomorphic function}.
\end{Remark}
\begin{Remark}
Let $(M,\mathcal{C})$ be a Lorentz surface and $(x,y),\ ({\tilde
x},{\tilde y})$ null coordinate systems. Then these two coordinate
systems are related by
$$
\frac{\partial \tilde x}{\partial y}=0,\ \ \frac{\partial \tilde
y}{\partial x}=0.
$$
Namely ${\tilde x}$ and ${\tilde y}$ depends only on $x$ and $y$
respectively.
\end{Remark}

\smallskip

On timelike surfaces of constant mean curvature $H\geq 0$, a
special (local) coordinate system is available (\cite{FI},
\cite{Hu}, \cite{Mil}, \cite{We}).

\begin{Proposition}
Let $\varphi:M\to {\mathbb E}^3_1$ be a timelike surface of
constant mean curvature $H\not=0$. Assume that $(M,\varphi)$ has
real distinct principal curvatures. Then there exists a local
coordinate system $(x,y)$ such that
$$
I=e^\omega \>{\mathrm d}x {\mathrm d}y ,\ \ II=\frac{H}{2} \left
\{ {\mathrm d}x^2+2e^\omega {\mathrm d}x{\mathrm d}y+{\mathrm
d}y^2 \right \}. \eqno(1.7)
$$
With respect to this local coordinate system, the Gauss-Codazzi
equation become
$$
\omega_{xy}+H^2\sinh \omega=0. \eqno{(\mathrm{shG})}
$$
\end{Proposition}
The partial differential equation (shG) is called the {\it
hyperbolic} {\it sinh-Gordon equation} or {\it affine Toda field
equation} of type $A^{(1)}_1$.

\begin{Remark}
Let $(x,y)$ be the local coordinate system in the preceeding
Proposition. Introduce a local coordinate system $(u,v)$ by
$x=u+v,\ y=-u+v$. Then $I$ and $II$ are represented as
$$
I=e^\omega (-du^2+dv^2),\ \ II=2He^{\frac{\omega}{2}} \left(
-\sinh \frac{\omega}{2}du^2+ \cosh \frac{\omega}{2}dv^2 \right).
$$
The local coordinate system $(u,v)$ is (Lorentz) isothermal and a
curvature-line coordinate system. Such a coordinate system $(u,v)$
is called an {\it isothermic coordinate system}.

Timelike CMC surfaces with real distinct principal curvatures are
called {\it isothermic timelike CMC surfaces}. Note that an
isothermic coordinate system is characterized as a local null
coordinate system $(x,y)$ such that $Q=R\not=0$. See \cite{FI}.

\end{Remark}

A very different situation is discussed in the
\begin{Proposition}
Let $\varphi:M\to {\mathbb E}^3_1$ be a timelike surface of
constant mean curvature $H\not=0$. Assume that $(M,\varphi)$ has
imaginary principal curvatures. Then there exists a local
coordinate system $(x,y)$ such that
$$
I=e^\omega {\mathrm d}x {\mathrm d}y ,\ \ II=\frac{H}{2} \left \{
{\mathrm d}x^2+2e^\omega {\mathrm d}x{\mathrm d}y -{\mathrm d}y^2
\right \}. \eqno(1.8)
$$
With respect to this local coordinate system, the Gauss-Codazzi
equations become
$$
\omega_{xy}+H^2\cosh \omega=0. \eqno{(\mathrm{chG})}
$$
\end{Proposition}
The partial differential equation (chG) is called the {\it
hyperbolic} {\it cosh-Gordon equation}.

The local coordinate system $(x,y)$ is called an {\it
anti-isothermic coordinate system}. The anti-isothermic coordinate
system is characterized as a local null coordinate system such
that $Q=-R\not=0$.

Timelike CMC surfaces with imaginary principal curvatures are
called {\it anti-isothermic timelike CMC surfaces}. The notion of
``anti-isothermic coordinate" has been introduced by \cite{FI}
(Definition 4.15).

\begin{Remark}
Let $\varphi:M\to {\mathbb E}^3_1$ be a timelike surface of
constant mean curvature $H$. Assume that $(M,\varphi)$ has two
equal and real principal curvatures. Then the Gauss equation of
$(M,\varphi)$ becomes the {\it Liouville equation}:
$$
\omega_{xy}+\frac{H^2}{2}e^\omega=0. \eqno{(\mathrm{L})}
$$
\end{Remark}

\vspace{0.2cm}

\noindent {\large {\bf 1.3}}\hspace{0.15cm} Next, we shall define
the Gau{\ss} map of a timelike surface. Let $M$ be a timelike
surface and $N$ a unit normal vector field to $M$. The {\it
Gau{\ss} map} $\psi$ of $M$ is a smooth map of $M$ into $S^2_1$,
which assigns to each $p \in M$, the point $\psi(p) \in
\mathbb{E}^3_1$ obtained by parallel translation of the unit
normal vector $N_{p}$ of $M$ at $p$ to the origin of
$\mathbb{E}^3_1$.

The constancy of the mean curvature is characterized by the
harmonicity of the Gauss map ({\it cf.\/} \cite{Mil}).

\begin{Proposition}
The Gauss map of a timelike surface is harmonic if and only if the
mean curvature is constant.
\end{Proposition}

\begin{Remark}

Let $(M_1,g_1)$ and $(M_1,g_2)$ be (semi-) Riemannian manifolds
and $\psi:M_1 \to M_2$ be a smooth map. Then $\psi$ is said to be
a {\it harmonic map} if its tension field $\tau(\psi)$:
$$
\tau(\psi):=\rm{tr}(\nabla {\mathrm d} \psi)
$$
vanishes. In case $\dim M_1=2$, the harmonicity of a smooth map
$\psi$ is invariant under the conformal transformation of
$(M_1,g_1)$. In particular, when $(M_1,g_1)$ is a Lorentzian
$2$-manifold, a harmonic map $\psi$ is often called a {\it
Lorentzian}) {\it harmonic map}.
\end{Remark}

Note that the constancy of the Gau{\ss}ian curvature is
characterized by the following.
\begin{Proposition}
Let $M$ be a timelike surface. Assume that the Gau{\ss}ian
curvature has a constant sign on $M$. Then the second fundamental
form $II$ gives $M$ another {\rm (}semi-{\rm )} Riemannian metric.
With respect to the conformal structure determined by $II$, the
Gau{\ss} map of $M$ is harmonic if and only if $K$ is constant.
\end{Proposition}

The Ad-action of $G$ on $S^2_1$ is transitive and isometric. The
isotropy subgroup $K$ of $G$ at ${\mathbf k}'$ is
$$
K=\left \{ u_0\,{\mathbf 1}+u_3\,{\mathbf k}'\ | \ u^2_0-u^2_3=-1
\right \}.
$$
The isotropy subgroup $K$ is isomorphic to the multiplicative
group $\mathbb{R}^*$. \par The natural projection $\pi :G=H^3_1
\rightarrow S^2_1$, given by $\pi(g)=\mathrm{Ad}(g){\mathbf k}'$,
$g \in  G$, defines a principal $\mathbb{R}^*$-bundle $H^3_1$ over
$S^2_1$. The fibering $\pi:H^3_1 \to S^2_1$ is called the {\it
Hopf-fibering of \/ }$S^2_1$.

The Lie algebra $\mathfrak{k}$ of $K$ is given by ${\mathfrak
k}={\mathbb R}\> {\mathbf k}'$. The tangent space of $S^2_1$ at
the origin ${\mathbf k}'$ is ${\mathfrak m}={\mathbb R}{\mathbf i}
\oplus {\mathbb R}{\mathbf j}'$. Let $\sigma$ be the involution of
${\mathfrak g}$ defined by $\sigma=\mathrm{Ad}({\mathbf
k}^{\prime})= \Pi_{\mathfrak k}-\Pi_{\mathfrak m}$, where
$\Pi_{\mathfrak k}$ and $\Pi_{\mathfrak m}$ are the projections
from ${\mathfrak g}$ onto ${\mathfrak k}$ and ${\mathfrak m}$
respectively. Then the pair $({\mathfrak g},\sigma)$ is a
symmetric Lie algebra data associated with  the Lorentzian
symmetric space $S^2_1=G/K$.

\section{Loop groups}
\noindent {\large {\bf 2.1}} \hspace{0.15mm} To study timelike CMC
surfaces in the spirit of \cite{DPW}, we need to introduce some
notation involving loop groups. Let us denote the polynomial loop
algebra of ${\mathfrak g}={\mathfrak s}{\mathfrak l}(2;{\mathbb
R})$ by $\Lambda_{\mathrm{pol}}\> \mathfrak{g}$:
$$
\Lambda_\mathrm{pol} \> \mathfrak{g} = \left \{
\xi(\lambda)=\sum_{\mathrm{finite}} \xi_{j}\lambda^{j}: S^1
\rightarrow \> {\mathfrak g} \right \}, \eqno(2.1)
$$

\noindent where  $S^1$ denotes the unit circle in $\mathbb{C}$.

Let $\sigma$ be the involution of $\mathfrak{g}$ corresponding to
the Lorentzian symmetric space $S^2_1=G/K$ defined above. Then the
{\it polynomial twisted loop algebra} of $\mathfrak{g}$ is defined
by
$$
\Lambda_\mathrm{pol}\> {\mathfrak g}_\sigma= \left \{ \xi(\lambda)
\in \Lambda_\mathrm{pol}{\mathfrak g} \ | \ \sigma( \xi(\lambda))
=\xi(-\lambda) \right \}. \eqno(2.2)
$$

For the purposes of this paper we need a certain Banach-completion
of $\Lambda_{\mathrm{pol}}\> {\mathfrak g}_\sigma$. To this end we
introduce the norm $|\cdot|_1$ for $\mathfrak{g}$:
$$
|A|_1:=\max_{j}\ \{\sum_{i=1}^{2} |a_{ij}|\ \}, \ \ A=(a_{ij}) \in
{\mathfrak g}.
$$
We extend this norm to the polynomial loop algebra in the
following way:

$$
|| \xi ||=\sum |\xi_j|_{1},\ \ \xi=\sum_{j} \xi_{j}\lambda^j \in
\Lambda_{\mathrm{pol}}\> {\mathfrak g}.
$$

Denote the completion of $\Lambda_{\mathrm{pol}}\>\mathfrak{g}$
and $\Lambda_{\mathrm{pol}}\> {\mathfrak g}_\sigma$ with respect
to the norm $||\cdot||$ by $\Lambda {\mathfrak g}$ and $\Lambda
{\mathfrak g}_\sigma$ respectively. Then the Lie algebras $\Lambda
{\mathfrak g}$ and $\Lambda {\mathfrak g}_\sigma$ are Banach Lie
algebras. ({\it cf.\/} Proposition 4.2.1 in \cite{T}.) Actually,
these are Banach algebras of continuous functions on $S^1$.
Moreover since the involution $\sigma$ is inner, these Banach Lie
algebras are isomorphic to each other.

Next we introduce the following Lie subalgebras of $\Lambda
\mathfrak{g}$:

$$
\Lambda^{+}{\mathfrak g} =\left \{ \xi(\lambda)=\sum_{j\geq 0}
\xi_{j}\lambda^{j} \in \Lambda {\mathfrak g} \right \},\ \
\Lambda^{-}{\mathfrak g} =\left \{ \xi(\lambda)=\sum_{j\leq 0}
\xi_{j}\lambda^{j} \in \Lambda {\mathfrak g} \right \}, \eqno(2.3)
$$
$$
\Lambda^{+}_{*} {\mathfrak g} =\left \{ \xi(\lambda)=\sum_{j> 0}
\xi_{j}\lambda^{j} \in \Lambda {\mathfrak g} \right \},\ \
\Lambda^{-}_{*} {\mathfrak g} =\left \{ \xi(\lambda)=\sum_{j< 0}
\xi_{j}\lambda^{j} \in \Lambda {\mathfrak g} \right \}. \eqno(2.4)
$$
Then we have the following decompositions as direct sums of linear
spaces:

$$
\Lambda {\mathfrak g}= \Lambda^{+}_{*} {\mathfrak g} \oplus
\Lambda^{-} {\mathfrak g} =\Lambda^{-}_{*} {\mathfrak g} \oplus
\Lambda^{+}{\mathfrak g}. \eqno(2.5)
$$

Similarly, we introduce the following Lie subalgebras of the
twisted loop algebra:

$$
\Lambda^{+}{\mathfrak g}_\sigma =\left \{ \xi(\lambda)=\sum_{j\geq
0} \xi_{j}\lambda^{j} \in \Lambda {\mathfrak g}_\sigma \right
\},\> \Lambda^{-}{\mathfrak g}_\sigma =\left \{
\xi(\lambda)=\sum_{j\leq 0} \xi_{j}\lambda^{j} \in \Lambda
{\mathfrak g}_\sigma \right \}. \eqno(2.6)
$$

$$
\Lambda^{+}_{*} {\mathfrak g}_\sigma =\left \{
\xi(\lambda)=\sum_{j> 0} \xi_{j}\lambda^{j} \in \Lambda {\mathfrak
g}_\sigma \right \},\> \Lambda^{-}_{*} {\mathfrak g}_\sigma =\left
\{ \xi(\lambda)=\sum_{j< 0} \xi_{j}\lambda^{j} \in \Lambda
{\mathfrak g}_\sigma \right \}. \eqno(2.7)
$$

Then we have the following decompositions as direct sums of linear
spaces:

$$
\Lambda {\mathfrak g}_{\sigma}= \Lambda^{+}_{*} {\mathfrak
g}_{\sigma} \oplus \Lambda^{-} {\mathfrak g}_{\sigma}
=\Lambda^{-}_{*} {\mathfrak g}_{\sigma} \oplus
\Lambda^{+}{\mathfrak g}_{\sigma}. \eqno(2.8)
$$

It is not difficult to see that one can analogously define
connected Banach Lie groups: $\Lambda G$, $\Lambda^{\pm}G$ and
$\Lambda^{\pm}_{*}G$, whose Lie algebras are $\Lambda
\mathfrak{g}$, $\Lambda^{\pm}\mathfrak{g}$ and
$\Lambda^{\pm}_{*}\mathfrak{g}$ respectively.

\

For the twisted case we have the following Banach Lie groups:
$\Lambda G_\sigma$, $\Lambda^{\pm}G_\sigma$ and
$\Lambda^{\pm}_{*}G_\sigma$, whose Lie algebras are $\Lambda
\mathfrak{g}_\sigma$, $\Lambda^{\pm}\mathfrak{g}_\sigma$ and
$\Lambda^{\pm}_{*}\mathfrak{g}_\sigma$ respectively.

\vspace{0.2cm}

\noindent {\large {\bf 2.2}} \hspace{0.15mm} In this subsection we
recall the classically known Birkhoff decomposition theorem for
loop groups of {\it complex} special linear group ${\rm
SL}(2;{\mathbb C})$.

\

We use this result in order to prove that a similar factorization
holds for the loop groups ${\tilde \Lambda}G$, and ${\tilde
\Lambda}G_{\sigma}$ defined below.

\begin{Theorem}
{\rm (}Birkhoff decomposition of $\Lambda G^{\mathbb C}${\rm )}
$$
\Lambda G^{\mathbb C} = \bigsqcup_{w \in \mathscr{T}} \Lambda
^{-}G^{\mathbb C} \cdot w \cdot \Lambda  ^{+}G^{\mathbb C},
\eqno{(2.9)}
$$
Here $\mathscr{T}$ denotes the group of homomorphisms from $S^{1}$
into the subgroup of diagonal matrices of
$\mathrm{SL}(2;\mathbb{C})$, that is,
$$
\mathscr{T}= \left \{ \> \left(
\begin{array}{cc}
\lambda^{a} & 0 \\
0 & \lambda^{-a}
\end{array}
\right ) \ \Biggr \vert \ a>0 \right \}.
$$
Moreover, the  multiplication maps
$$\Lambda ^{-}_{*}G^{\mathbb C}
\times
 \Lambda  ^{+}G^{\mathbb C}
\to   \Lambda G^{\mathbb C},\ \ \Lambda^{+}_{*}G^{\mathbb C}
\times \Lambda   ^{-}G^{\mathbb C} \to   \Lambda G^{\mathbb C}.
\eqno(2.10)
$$
are diffeomorphisms onto the open dense subsets
$\mathscr{B}^{\circ}_{\Lambda}(-,+)$ and
$\mathscr{B}^{\circ}_{\Lambda}(+,-)$ of $\Lambda G^{\mathbb C}$,
called the big cells of $\Lambda  G^{\mathbb C}$. In particular if
$\gamma$ is an element of $\mathscr{B}_{\Lambda}
=\mathscr{B}^{\circ}_{\Lambda}(-,+) \cap
\mathscr{B}^{\circ}_{\Lambda}(+,-)$, then $\gamma$ has unique
decompositions{\rm :}
$$
\gamma=\gamma_{-}\cdot\ell_{+} =\gamma_{+}\cdot\ell_{-}, \ \
\gamma_{\pm}\in \Lambda^{\pm}_{*}G^{\mathbb C},\ \ell_{\pm}\in
\Lambda^{\pm}G^{\mathbb C}.
$$
Here the subgroups $ \Lambda^{\pm}_{*} G^{\mathbb C}$ are defined
by
$$
 \Lambda ^{-}_{*} G^{\mathbb C}=
\{ \ \gamma \in
 \Lambda  ^{-}G^{\mathbb C}
\ | \ \gamma(\lambda)=\mathbf{1} +\sum_{k\leq -1}\gamma_{k}
\>\lambda^k \ \},
$$
$$
\Lambda^{+}_{*} G^{\mathbb C}= \{ \ \gamma \in
 \Lambda^{+}G^{\mathbb C}
\ | \ \gamma(\lambda)=\mathbf{1} +\sum_{k\geq 1}\gamma_{k}
\>\lambda^k \ \}.
$$
\end{Theorem}

Next let ${\tilde \Lambda}G$ be the subset of $\Lambda G$ whose
elements, as maps defined on $S^1$, admit analytic continuations
to $\mathbb{C}^*$.

$$
{\tilde \Lambda}G= \left \{ \gamma \in \Lambda G \ | \ \gamma: S^1
\rightarrow G \ \mathrm{extends}\ \mathrm{analytically} \
\mathrm{to} \ {\mathbb C}^{*} \ \right \}. \eqno(2.11)
$$
Similarly we define
$$
{\tilde \Lambda}G_\sigma= \left \{ \gamma \in \Lambda G_\sigma\ |
\ \gamma \ \mathrm{extends}\ \mathrm{analytically} \ \mathrm{to} \
{\mathbb C}^{*} \ \right \}. \eqno(2.12)
$$

It is easy to check that ${\tilde \Lambda} G$ is a subgroup of
$\Lambda G$. Similarly, ${\tilde \Lambda} G_{\sigma}$ is a
subgroup of $\Lambda G_{\sigma}$.

For ${\tilde \Lambda} G$ we will use the topology induced  from ${
\Lambda} G$. Then we obtain the following theorem.

\begin{Theorem}
{\rm (}Birkhoff decomposition of ${\tilde \Lambda}G${\rm )}
$$
{\tilde \Lambda} G = \bigsqcup_{w \in \mathscr{T}} {\tilde
\Lambda}^{-}G \cdot w \cdot {\tilde \Lambda}^{+}G.
$$
Here $\mathscr{T}$ denotes the group of homomorphisms from $S^1$
into the subgroup of diagonal matrices of
$\mathrm{SL}(2;\mathbb{C})$. Moreover, the  multiplication maps
$$
{\tilde \Lambda}^{-}_{*}G \times {\tilde \Lambda}^{+}G \to {\tilde
\Lambda} G, \ \ \ {\tilde \Lambda}^{+}_{*}G \times {\tilde
\Lambda}^{-}G \to {\tilde \Lambda} G.
$$
are diffeomorphism onto the open dense subsets
$\mathscr{B}^{\circ}(-,+)$ and $\mathscr{B}^{\circ}(+,-)$ of
${\tilde \Lambda} G$, called the big cells of ${\tilde \Lambda}G$.
In particular if $\gamma$ is an element of
$\mathscr{B}^{\circ}=\mathscr{B}^{\circ}(-,+) \cap
\mathscr{B}^{\circ}(+,-)$, then $\gamma$ has unique
decompositions{\rm :}
$$
\gamma=\gamma_{-}\cdot\ell_{+} =\gamma_{+}\cdot\ell_{-}, \ \
\gamma_{\pm}\in {\tilde \Lambda}^{\pm}_{*}G,\ \ell_{\pm}\in
{\tilde \Lambda}^{\pm}G.
$$
\end{Theorem}

\begin{pf}
The idea of this proof is similar to the one presented in the
Appendix of \cite{T}.

Let $g \in {\tilde \Lambda} G$ with expansion $g(\lambda)=\sum
g_{j}\> \lambda^{j}$.

Note that the coefficients $g_{j}$ in the expansion of $g$ are
real.

Over the unit circle $S^1$, we obtain a decomposition $g=g_{-}
\cdot w \cdot g_{+}$ by the classical Birkhoff Decomposition
Theorem 2.2. It remains to show that actually every factor of $g$
defines an element in ${\tilde \Lambda} G$.

First we show all factors of $g$ are in $\Lambda G$. To see this
we introduce the automorphism $\kappa$ of $\Lambda G^{\mathbb{C}}$
which is defined by
$$
\kappa(\gamma)(\lambda):= \sum {\bar \gamma}_{j}\lambda^{j}
$$
for every
$$
\gamma(\lambda)=\sum \gamma_{j}\lambda^{j} \in \Lambda G^{\mathbb
C}.
$$
It is obvious to verify that $\kappa$ leaves $\Lambda^{+}
G^{\mathbb{C}}$ and $\Lambda^{-} G^{\mathbb{C}}$ invariant and
fixes all $w \in \mathscr{T}$. Moreover, $\Lambda G$ is the fixed
point set of $\kappa$. Thus $g = \kappa (g) = \kappa(g_{-}) \cdot
w \cdot \kappa(g_{+})$. The classical Birkhoff Decomposition
Theorem 2.2 now implies $\kappa(g_{-}) = g_{-} \>\cdot \>v_{-}$
and $\kappa(g_{+}) = v_{+} \> \cdot \>g_{+}$, where $v_{\pm} \in
\Lambda^{\pm}G^{\mathbb C}$. Moreover, we have $v_{-} \cdot w = w
\cdot v_{+}^{-1}$.

Applying $\kappa$ again to $\kappa(g_{-})=g_{-}\>v_{-}$ and taking
into account that $\kappa$ is an automorphism of order two, it
follows that
$$
g_{-}=\kappa (\kappa(g_{-}))= \kappa(g_{-}\>v_{-})= g_{-}\cdot
v_{-} \cdot \kappa(v_{-}).
$$
Thus we obtain $\kappa(v_{-})=(v_{-})^{-1}$.

Analogusly, from $\kappa(g_{+})=v_{+}\>g_{-}$, we obtain $\kappa
(v_{+}) =  (v_{+})^{-1}$.

If $w = \mathbf{1}$, then we are done. Assume now $w \ne
\mathbf{1}$. Then $v_{-}$ and $v_{+}$ are lower triangular.
Moreover, the equation relating $v_{-}$ and $v_{+}$ via $w$ shows
that the diagonal part of $v_{-}$ is independent of $\lambda$.
Thus $\kappa (v_{-}) = (v_{-})^{-1}$ shows that the diagonal
entries of $v_{-}$ have  modulus $1$. Let $d$ denote the diagonal
part of $v_{-}$ and let $\sqrt{d}$ denote its square root. Then we
set $g^{\prime}_{-}:=g_{-} \cdot (\sqrt{d})^{-1}$ and obtain
$\kappa(g^{\prime}_{-}) =g^{\prime}_{-} \cdot v^{\prime}_{-}$,
where $v^{\prime}_{-}$ is lower triangualr with $1$'s on the
diagonal. Writing now $g = g^{\prime}_{-} \cdot w \cdot
g^{\prime}_{+}$, where $g^{\prime}_{+}$ is defined by
$g^{\prime}_{+}= \sqrt{d} \cdot g_{+}$, then we see that the
corresponding $v'_{+}$ has $1$'s on the diagonal. Now we take the
square root of $v^{\prime}_{-}$ and multiply $g^{\prime}_{-}$ by
its inverse on the right, obtaining $g^{\prime \prime}_{-}$. A
straightforward computation shows that $\kappa$ fixes $g^{\prime
\prime}_{-}$, whence  $g^{\prime \prime}_{-}$ is in $\Lambda G$.
As a consequence, the corresponding $g^{\prime \prime}_{+}$ is
also in $\Lambda G$.

\vspace{0.2cm}

Next we need to show that $g_{-}$ and $g_{+}$ are in the `` $\sim$
''-group.

We know from the definition that $g_{-}$ has a holomorphic
extension to the exterior of the unit disk and is finite at
$\infty$. Thus $(g_{-})^{-1} \cdot g = w \cdot g_{+}$ has a
holomorphic extension to the exterior of the unit disk. Hence
$g_{+}$ also has a holomorphic extension to the exterior of the
unit disk. Altogether, $g_{+}$ has a holomorphic extension to
$\mathbb{C}^*$. This proves the first part of the Theorem.
\newline
\indent For the second part we note first that the big cells are
indeed dense, since the cosets involving $w \neq \mathbf{1}$ are
of nonzero codimension in in $\Lambda G$.

The rest of the claim follows by the fact the multiplication maps
are induced by the diffeomorphisms $ \Lambda ^{-}_{*}G^{\mathbb C}
\times\Lambda   ^{+}G^{\mathbb C} \to  \Lambda G^{\mathbb C}$ and
$ \Lambda ^{+}_{*}G^{\mathbb C} \times \Lambda  ^{-}G^{\mathbb C}
\to \Lambda  G^{\mathbb C}$ via the loop group correspondences
those we presented in this section. $\Box$
\end{pf}

Finally, we consider the twisted loop groups defined earlier. We
also define the twisted analytic loop groups derived from ${\tilde
\Lambda} G_{\sigma}$ in the obvious way. Then we have

\begin{Theorem} {\rm (}Birkhoff decomposition
of ${\tilde \Lambda}G_\sigma$ {\rm )} Theorem {\rm 2.2} also holds
for the twisted groups.
\end{Theorem}

\begin{pf} We note that the twisting involution
$\sigma$ is given by an inner automorphism of $G$. Therefore, the
twisted loop group and the untwisted loop group are isomorphic.
Actually, the isomorphism from the untwisted loop group to the
twisted loop group is given easily: powers $\lambda^k$ on the
diagonal are doubled, powers $\lambda^k$ in the $(1,2)$-position
are replaced by $\lambda^{2k+1}$ and in the $(2,1)$-position they
are replaced by $\lambda^{2k-1}$. The claim now follows. $\Box$
\end{pf}
\begin{Remark} The proof above actually also shows the
Birkhoff Decomposition Theorem for $\Lambda G_{\sigma}$ .
\end{Remark}

\vspace{0.2cm}

\noindent {\large {\bf 2.3}} \hspace{0.15mm} We have the following
fundamental decomposition theorem:

\begin{Theorem}
{\rm (}Iwasawa decomposition of ${ \Lambda}G_\sigma \times {
\Lambda} G_\sigma${\rm )}

Let $\Delta ( { \Lambda} G_\sigma \times { \Lambda} G_\sigma)$
denote the diagonal subgroup of $ { \Lambda} G_\sigma \times {
\Lambda} G_\sigma$. Then we have

$$
{ \Lambda} G_{\sigma} \times { \Lambda} G_{\sigma} = \bigsqcup
\Delta( { \Lambda} G_{\sigma} \times { \Lambda} G_{\sigma}) \cdot
({\mathbf 1},w) \cdot ({ \Lambda}^{-} G_{\sigma} \times {
\Lambda}^{+} G_{\sigma}),
$$
where $w \in \mathscr{T}$ is as in Theorem {\rm 2.2}.

Moreover, the multiplication maps
$$
\Delta ({ \Lambda} G_{\sigma} \times { \Lambda} G_{\sigma}) \times
({ \Lambda}^{-}_{*}G_{\sigma} \times { \Lambda}^{+}G_{\sigma}) \to
{ \Lambda}G_{\sigma} \times { \Lambda}G_{\sigma},
$$
$$
\Delta ({ \Lambda} G_{\sigma} \times { \Lambda} G_{\sigma}) \times
({ \Lambda}^{+}_{*}G_{\sigma}
 \times
{ \Lambda}^{-}G_{\sigma}) \to { \Lambda}G_{\sigma} \times {
\Lambda}G_{\sigma}
$$
are diffeomorphisms onto the open dense subsets
$\mathscr{I}_{\Lambda}(+,-)$ and $\mathscr{I}_{\Lambda}(-,+)$ of
${ \Lambda}G_{\sigma} \times { \Lambda}G_{\sigma}$-- called the
big cells of ${ \Lambda}G_{\sigma} \times { \Lambda}G_{\sigma}$ .
\end{Theorem}

\begin{pf}
Take $(g,h) \in { \Lambda}G_{\sigma}\times { \Lambda}G_{\sigma}$.
Decompose $g^{-1}\cdot h$ acoording to the Birkhoff decomposition
of ${ \Lambda}G_{\sigma}$ (Theorem 2.2):
$$
g^{-1}\cdot h=u_{-}\> w\> u_{+},\ \ u_{\pm} \in {
\Lambda}G_{\sigma},\ w \in \mathscr{T}.
$$
It is easy to verify that the splitting
$$
(g,h)=(gu_{-}, gu_{-})({\mathbf 1},w) (\>(u_{-})^{-1},u_{+})
$$
gives the Iwasawa decomposition of $(g,h)$ in the {\it
un\/}twisted loop group ${ \Lambda}G \times { \Lambda}G$. ({\it
cf.\/} Theorem 4.1 in \cite{BG}.) Since the factors $gu_{-}$ and
$u_{\pm}$ are $\sigma$-twisted, this splitting is the required
(Iwasawa) splitting in the {\it twisted} loop group ${
\Lambda}G_{\sigma} \times { \Lambda}G_{\sigma}$.

For the second claim we note that our definitions imply
$$
\Delta ({ \Lambda} G_{\sigma} \times { \Lambda} G_{\sigma}) \
\bigcap \ ({ \Lambda}^{-}_{*}G_{\sigma}
 \times
{ \Lambda}^{+}G_{\sigma}) =\{\mathbf{1} \}.
$$

Thus the splitting is unique, whence the map is a bijection onto
its image. The proof that it is  a diffeomorphism is almost
verbatim the same as for the Birkhoff decomposition. In fact, the
proof follows as in \cite{DGS}.

It remains to show that the big cell is dense. But if $(g,h)$ is
given, then $g^{-1} \cdot h$ is in $\Lambda G_{\sigma} \cong
\Lambda G$ and from the Birkhoff decomposition Theorem 2.3, we
know that in every neighbourhood of this element there is an
element in the big cell. Therefore, in every neighbourhood of $g$
and $h$ there exists some $g'$ and $h'$ such that $(g')^{-1} \cdot
h'$ is in the big cell of $\Lambda G_{\sigma}$. But the proof
above shows that then $(g',h')$ is in the big cell of $\Lambda
G_{\sigma} \times \Lambda G_{\sigma}$. $\Box$
\end{pf}

\vspace{0.2cm}

\noindent {\large {\bf 2.4}} \hspace{0.15mm} For twisted loop
groups of elements with analytic extension, we have the following
decomposition theorem:

\begin{Theorem}
{\rm (}Iwasawa decomposition of ${\tilde \Lambda}G_\sigma \times
{\tilde \Lambda} G_\sigma${\rm )}

Let $\Delta ( {\tilde \Lambda} G_\sigma \times {\tilde \Lambda}
G_\sigma)$ denote the diagonal subgroup of $ {\tilde \Lambda}
G_\sigma \times {\tilde \Lambda} G_\sigma$. Then we have

$$
{\tilde \Lambda} G_{\sigma} \times {\tilde \Lambda} G_{\sigma} =
\bigsqcup \Delta( {\tilde \Lambda} G_{\sigma} \times {\tilde
\Lambda} G_{\sigma}) \cdot ({\mathbf 1},w) \cdot ({\tilde
\Lambda}^{-} G_{\sigma} \times {\tilde \Lambda}^{+} G_{\sigma}),
$$
where $w \in \mathscr{T}$ is as in Theorem {\rm 2.2}.

Moreover, the multiplication maps
$$
\Delta ({\tilde \Lambda} G_{\sigma} \times {\tilde \Lambda}
G_{\sigma}) \times ({\tilde \Lambda}^{-}_{*}G_{\sigma} \times
{\tilde \Lambda}^{+}G_{\sigma}) \to {\tilde \Lambda}G_{\sigma}
\times {\tilde \Lambda}G_{\sigma},
$$
$$
\Delta ({\tilde \Lambda} G_{\sigma} \times {\tilde \Lambda}
G_{\sigma}) \times ({ \Lambda}^{+}_{*}G_{\sigma}
 \times
{\tilde \Lambda}^{-}G_{\sigma}) \to {\tilde \Lambda}G_{\sigma}
\times {\tilde \Lambda}G_{\sigma}
$$
are diffeomorphisms onto the open dense subsets $\mathscr{I}(+,-)$
and $\mathscr{I}(-,+)$ of ${\tilde \Lambda}G_{\sigma} \times
{\tilde \Lambda}G_{\sigma}$-- called the big cells of ${\tilde
\Lambda}G_{\sigma} \times {\tilde \Lambda}G_{\sigma}$ .
\end{Theorem}
\begin{pf}
Take $(g,h) \in {\tilde \Lambda}G_{\sigma}\times {\tilde
\Lambda}G_{\sigma}$. Decompose $g^{-1}\cdot h$ acoording to the
Birkhoff decomposition of ${\tilde \Lambda}G_{\sigma}$ (Theorem
2.2):
$$
g^{-1}\cdot h=u_{-}\> w\> u_{+},\ \ u_{\pm} \in {\tilde
\Lambda}G_{\sigma},\ w \in \mathscr{T}.
$$
As we showed in the proof of Theorem 2.4, the splitting
$$
(g,h)=(gu_{-}, gu_{-})({\mathbf 1},w) (\>(u_{-})^{-1},u_{+})
$$
gives the Iwasawa decomposition of $(g,h)$ in the {\it twisted}
loop group ${\tilde \Lambda}G_{\sigma} \times {\tilde
\Lambda}G_{\sigma}$.

We need to check that $u_{-}$ and $u_{+}$ are in the
``$\sim$"-group. By definition, $u_{-}$ has a holomorphic
extension to the exterior of the unit disk and is finite at
$\infty$. Thus $u_{+}=(gu_{-}w)^{-1}h$ has  a holomorphic
extension to the exterior of the unit disk. Thus $u_{+}$ has
holomorphic extension to $\mathbb{C}^{*}$. Similarly $u_{-}$ has
also a holomorphic extension to $\mathbb{C}^{*}$. This proves the
claim. The remaining assertions follow from the previous Theorem,
since we use the induced topology. $\Box$

\end{pf}

\begin{Remark}
P.~Kellersch \cite{K} generalized the classical Iwasawa
decomposition for untwisted loop groups of compact simple Lie
groups to those for loop groups of general simple Lie groups.
Moreover V.~Balan and the first named author \cite{BD} generalized
the splitting theorem due to Kellersch to those for general Lie
groups.

\end{Remark}

\begin{Remark}
I.~T.~Gohberg and his collaborators investigated splittings for
matrix valued functions over (separate) contours and more general
Banach algebras of functions. For instance, the Birkhoff splitting
for matrix valued functions with coefficients in the {\it Wiener
algebra}:
$$
\mathscr{A}= \left \{ f(\lambda)=\sum f_j\> \lambda^{j}:S^1 \to
\mathbb{C}\ \biggr | \ \sum |f_{j}|<\infty \right \}
$$
was proven by Gohberg in \cite{Go}. But also
 splittings for matrix
valued functions with coefficients in the Wiener algebra on the
{\it real line} were obtained \cite{CG}. For more  information, we
refer to \cite{Go} and \cite{CG} references therein.

For our geometric purposes--a Weierstra{\ss} type representation
for timelike surfaces in Minkowskispace --, we need {\it real}
Banach algebras of functions. In fact, we  need actually a {\it
real} loop parameter $\lambda$.

Fortunately all the loop group elements occurring in our geometric
context (extended framings etc.) have analytic extensions to
$\mathbb{C}^{*}$, i.e. these geometric loop group elements are all
contained in ${\tilde \Lambda} G$. For this reason we have
presented in this section the Birkhoff and Iwasawa decomposition
theorems for the twisted analytic loop group ${\tilde \Lambda}
G_{\sigma}$.
\end{Remark}

\section{Harmonic maps into $S^2_1$}

{\large {\bf 3.1}}\hspace{0.15cm} In this section we shall derive
a correspondence between harmonic maps from a simply connected
Lorentz surface $\mathbb{D}$ to $S^2_1$ and certain kind of flat
connections. (so-called {\it zero curvature representation}).
\par
We note that since the harmonic map equation is a local condition,
it suffices to consider harmonic maps from
 simply connected Lorentz surfaces
into $S^2_1$.

For the rest of this paper ${\mathbb D}$ will always denote a
simply connected region of the Minkowski plane
$(\mathbb{R}^2(x,y), \mathrm{d}x\mathrm{d}y)$ containing the
origin.

The following result is the starting point of our approach.

\begin{Proposition}
A smooth map $\psi : {\mathbb D} \rightarrow S^2_1 \subset
{\mathbb E}^3_1$ is harmonic if and only if
$$
\frac{{\partial}^2 \psi}{\partial x \partial y}= \rho \>\psi
\eqno(3.1)
$$
for some function $\rho$ on ${\mathbb D}$.
\end{Proposition}

\vspace{0.2cm}

\noindent {\large {\bf 3.2}}\hspace{0.15cm} Let $\psi:{\mathbb
D}\to S^2_1$ be a smooth map and $\pi:G \to S^2_1$ the Hopf
fibration as before. Since $\mathbb{D}$ is simply connected,
$\psi$ has a smooth lift $\Psi:\mathbb{D}\to G$ unique up to the
right $K$-action. Such a lift $\Phi$ is called a {\it framing} of
$\psi$. Note that a framing $\Psi$ is related to $\psi$ by
$$
\psi=\mathrm{Ad}(\Psi)\mathbf{k}^\prime. \eqno(3.2)
$$
Since $\mathbb{D}$ is simply connected, the pull-back bundle
$\psi^{*}\ G$ is necessarily a trivial bundle ${\mathbb D} \times
K$. So the group ${\mathcal G}$ of gauge transformations of
$\psi^{*}\ G$ is identified with ${\mathrm C}^{\infty}({\mathbb
D},K)$. \par We wish to describe the harmonicity of $\psi$ in
terms of a framing. Let $\mu_G$ be the Maurer-Cartan form of $G$.
It is well-known that $\mu_G$ satisfies the {\it Maurer-Cartan
equation}:
$$
{\mathrm d} \mu_G+\frac{1}{2}[\mu_G \wedge \mu_G]=0.
$$
The pulled back $1$-form $\alpha=\Psi^* \mu_G=\Psi^{-1}
\mathrm{d}\Psi$ of $\mu_G$ by $\Psi$ then satisfies
$$
{\mathrm d}\alpha+\frac{1}{2}[\alpha \wedge \alpha]=0. \eqno(3.3)
$$

The identity $(3.3)$ is equivalent with the integrability
condition for the existence of a smooth map $\Psi : {\mathbb D}
\rightarrow G$ such that $\alpha=\Psi^*\mu_G$. (Frobenius
theorem). By definition, $\alpha$ is a ${\mathfrak g}$-valued
$1$-form on $\mathbb{D}$. The ${\mathfrak g}$-valued $1$-form
$\alpha$ has a type decomposition along the decomposition
${\mathfrak g}={\mathfrak k} \oplus {\mathfrak m}$;
$$
\alpha=\alpha_0+\alpha_1. \eqno(3.4)
$$

Here $\alpha_0$ and  $\alpha_1$ denote the ${\mathfrak k}$-valued
part and ${\mathfrak m}$-valued part respectively. Write
$$
\alpha_{0}=\alpha_{\mathfrak k}^{\prime}\>{\mathrm d}x+
\alpha_{\mathfrak k}^{\prime \prime}\>{\mathrm d}y,\ \
\alpha_{1}=\alpha_{\mathfrak m}^{\prime}\>{\mathrm d}x+
\alpha_{\mathfrak m}^{\prime \prime}\>{\mathrm d}y. \eqno(3.5)
$$
Then $\alpha_0$ and $\alpha_1$ are decomposed with respect to the
conformal structure of $\mathbb{D}$ as follows:
$$
\alpha_0=\alpha_0^{\prime}+ \alpha_0^{\prime \prime}, \
\alpha_1=\alpha_1^{\prime}+ \alpha_1^{\prime \prime}, \eqno(3.6)
$$
$$
\alpha_0^{\prime}= \alpha_{\mathfrak k}^{\prime}\> {\mathrm d}x,\
\ \alpha_0^{\prime \prime}= \alpha_{\mathfrak k}^{\prime \prime}\>
{\mathrm d}y, \ \ \alpha_1^{\prime}=\alpha_{\mathfrak
m}^{\prime}\> {\mathrm d}x,\ \ \alpha_1^{\prime \prime}=
\alpha_{\mathfrak m}^{\prime \prime}\> {\mathrm d}y. \eqno(3.7)
$$

Define $\alpha^{\prime}$ and $\alpha^{\prime \prime}$ by
$$
\alpha^{\prime}:=\alpha^{\prime}_{0}+\alpha^{\prime}_{1}, \
\alpha^{\prime \prime}:= \alpha^{\prime \prime}_{0}+\alpha^{\prime
\prime}_{1}. \eqno(3.8)
$$
The 1-forms  $\alpha^{\prime}$ and $\alpha^{\prime \prime}$ are
called the (1,0)-{\it part} and (0,1)-{\it part} of $\alpha$
respectively.

\vspace{0.2cm}

With respect to the conformal structure of ${\mathbb D}$, we
decompose the exterior differential operator ${\mathrm d}$ ({\it
cf.\/} $(1.6)$):
$$
{\mathrm d}={\mathrm d}^{\prime}+ {\mathrm d}^{\prime \prime},\ \
{\mathrm d}^{\prime}= \frac{\partial}{\partial x}{\mathrm d}x,\ \
{\mathrm d}^{\prime \prime} =\frac{\partial}{\partial y} {\mathrm
d}y. \eqno(3.9)
$$
Then we have
$$
\alpha^{\prime}=\Psi^{-1}{\mathrm d}^{\prime}\Psi,\  \
\alpha^{\prime \prime}=\Psi^{-1}{\mathrm d}^{\prime \prime}\Psi.
\eqno(3.10)
$$

By the usual computations we obtain the following ({\it cf.}
\cite{GO}, \cite{MS})
\begin{Proposition}
Let $\psi : {\mathbb D}\rightarrow S^2_1$ be a smooth map with
framing $\Psi$. Then $\psi$ is harmonic if and only if
$$
{\mathrm d}(*\alpha_1)+[\alpha_0 \wedge * \alpha_1]=0 \eqno(3.11)
$$
for $\alpha=\Psi^{-1} {\mathrm d}\Psi.$\par Here $*$ denotes the
Hodge star operator acting on {\rm (}${\mathfrak g}$-valued {\rm
)} one forms on ${\mathbb D}$ defined by
$$
*\ {\mathrm d}x={\mathrm d}x,\ \ *\ {\mathrm d}y=-{\mathrm d}y.
\eqno(3.12)
$$
\end{Proposition}

\noindent {\large {\bf 3.3}}\hspace{0.15cm} Let $A^{1}({\mathbb
D};{\mathfrak g})$ be the space of all ${\mathfrak g}$-valued
one-forms and ${\mathcal A}$ the affine space of all connection
1-forms on the  product bundle ${\mathbb D}\times G$. The space
${\mathcal A}$ is an affine space associated to the linear space
$A^{1}({\mathbb D};{\mathfrak g})$. We shall choose the trivial
flat connection as the origin of ${\mathcal A}$, then the space
${\mathcal A}$ is identified with $A^{1}({\mathbb D};{\mathfrak
g})$. Hereafter we shall identify ${\mathcal A}$ with
$A^{1}({\mathbb D};{\mathfrak g})$ in this way.
\begin{Definition}
A connection $\alpha \in {\mathcal A}=A^{1}({\mathbb D};{\mathfrak
g})$ is {\it admissible} provided that
$$
{\mathrm d}\alpha+\frac{1}{2}[\alpha \wedge \alpha]=0,\ \ {\mathrm
d}(*\alpha_1)+[\alpha_0 \wedge * \alpha_1]=0. \eqno(3.13)
$$
\end{Definition}

The space of all admissible connections on ${\mathbb D}$ is
denoted by ${\mathcal A}^{*}$. The Frobenius' theorem implies the
following:
\begin{Lemma}
Let $\alpha$ be an admissible connection. Then there exists a
harmonic map $\psi:{\mathbb D}\to S^2_1$ such that for any framing
$\Psi$ of $\psi$ we have $\Psi^{-1}{\mathrm d}\Psi=\alpha$.

\end{Lemma}

Recall that framings are unique up to multiplication on the left
by  matrices independent of $x$ and $y$. We will remove this
freedom in the loop group formalism discussed below.

\vspace{0.2cm}

\noindent {\large {\bf 3.4}}\hspace{0.15cm} Furthermore for an
admissible connection $\alpha$, we define the the scalar field
${\mathscr S}(\alpha)$ by
$$
{\mathscr S}(\alpha)(x,y):= \langle \alpha_{\mathfrak
m}^{\prime}(x,y),\ \alpha_{\mathfrak m}^{\prime}(x,y) \rangle
\cdot \langle \alpha_{\mathfrak m}^{\prime \prime}(x,y),\
\alpha_{\mathfrak m}^{\prime \prime}(x,y) \rangle. \eqno(3.14)
$$

Let us define the linear subspaces ${\mathcal A}_{+},\ {\mathcal
A}_{0},\ {\mathcal A}_{-}$ of ${\mathcal A}^{*}$ by
$$
{\mathcal A}_{+}:= \{ \alpha \in {\mathcal A}^{*} \ \vert \
{\mathscr S}(\alpha)>0 \},\ \ {\mathcal A}_{-}:= \{ \alpha \in
{\mathcal A}^{*} \ \vert \ {\mathscr S}(\alpha)<0 \}, \eqno(3.15)
$$
$$
{\mathcal A}_{0}:= \{ \alpha \in {\mathcal A}^{*} \ \vert \
{\mathscr S}(\alpha)=0 \}.
$$
It is easily checked that all the spaces; ${\mathcal A}^{*}$,
${\mathcal A}_{\pm}$ and ${\mathcal A}_{0}$ are invariant under
the action of the gauge group
 as well as under conformal changes of
${\mathbb D}$. Note that the gauge  group ${\mathcal G}$ acts on
${\mathcal A}$ as follows:
$$
g^{*}\alpha=g^{-1}\ dg+{\mathrm A}{\mathrm d}(g^{-1})\alpha
\eqno(3.16)
$$
for $g \in {\mathcal G},$ $\alpha \in {\mathcal A}$.

\begin{Proposition}
Let ${\mathcal H}^{*}$ be the space of all harmonic maps from
${\mathbb D}$ to $S^2_1$. Define the linear subspaces ${\mathcal
H}_{+}$, ${\mathcal H}_{-}$ and ${\mathcal H}_{0}$ of ${\mathcal
H}^{*}$ by
$$
{\mathcal H}_{+}= \{ \psi \in {\mathcal H}^{*}\ | \ \langle
\psi_x,\> \psi_x \rangle \ \langle \psi_y,\> \psi_y \rangle >0 \
\mathrm{on}\  {\mathbb D} \}, \eqno(3.17)
$$
$$
{\mathcal H}_{-}= \{ \psi \in {\mathcal H}^{*}\ | \ \langle
\psi_x,\> \psi_x \rangle \ \langle \psi_y,\> \psi_y \rangle <0 \
\mathrm{on}\  {\mathbb D} \}, \eqno(3.18)
$$
$$
{\mathcal H}_{0}= \{ \psi \in {\mathcal H}^{*}\ | \ \langle
\psi_x,\> \psi_x \rangle \ \langle \psi_y,\> \psi_y \rangle =0 \
\mathrm{on}\  {\mathbb D}\}. \eqno(3.19)
$$
Then there are the following bijective correspondences{\rm :}
$$
{\mathcal H}^{*}\longleftrightarrow {\mathcal A}^{*}/{\mathcal G},
\ \ {\mathcal H}_{\pm}\longleftrightarrow {\mathcal
A}_{\pm}/{\mathcal G},\ \ {\mathcal H}_{0}\longleftrightarrow
{\mathcal A}_{0}/{\mathcal G}. \eqno(3.20)
$$
Here ${\mathcal G}$ denotes the gauge transformation group of the
bundle $\psi^{*}G={\mathbb D} \times K$. These correspondences are
described by $\alpha=\Psi^{-1}{\mathrm d}\Psi$ via a framing
$\Psi$ of $\psi$.
\end{Proposition}
\begin{pf}
It suffices to show that a harmonic map $\psi \in {\mathcal
H}_{\pm}$ [resp. ${\mathcal H}_{0}$] corresponds to a gauge class
of an admissible connection with $\pm{\mathscr S}(\alpha)>0$
[resp. ${\mathscr S}(\alpha)=0$]. However this is apparent from
the relations:
$$
\frac{\partial}{\partial x}\psi={\mathrm A}{\mathrm d}(\Psi)
[\alpha^{\prime}_{\mathfrak m},\ {\mathbf k}^\prime], \ \ \
\frac{\partial}{\partial y}\psi={\mathrm A}{\mathrm d}(\Psi)
[\alpha^{\prime \prime}_{\mathfrak m},\ {\mathbf k}^\prime].
$$
$\Box$
\end{pf}

\begin{Remark}
Here is a differential geometric interpretation of the subspaces
${\mathcal H}_{\pm},\ {\mathcal H}_{0}$.

We know that to any harmonic map $\psi \in {\mathcal H}_{\pm}$,
there exists a timelike immersion $\varphi$. (See Section 4.) The
condition $\psi \in {\mathcal H}_{+}$ [resp. $\psi \in {\mathcal
H}_{-}]$ is equivalent to the positivity [resp. negativity] of the
discriminant for the characteristic equation of the shape operator
of $\varphi$. Of course the condition $\psi \in {\mathcal H}_{0}$
corresponds to the property ``$\varphi$ has real repeated
principal curvatures".
\end{Remark}

\vspace{0.2cm}

\noindent {\large {\bf 3.5}}\hspace{0.15cm} To close this section,
we introduce the so-called {\it spectral parameter}.
\begin{Definition}
Let $\alpha \in {\mathcal A}$ be a connection. {\it A loop}
$\alpha_{\lambda}$ {\it of connections through} $\alpha$ is
defined by the following rule:
$$
\alpha_{\lambda}=\alpha_0+\lambda \alpha_{1}^{\prime}+
\lambda^{-1}\alpha_{1}^{\prime \prime},\ \ \lambda \in {\mathbb
R}^{+} \eqno(3.21)
$$
\end{Definition}
Note that $\mathscr{S}(\alpha_\lambda)\equiv \mathscr{S}(\alpha)$.
\newline
\indent The following observation is fundamental for our approach.

\begin{Proposition}
A connection $\alpha \in {\mathcal A}$ is admissible if and only
if the loop $\alpha_{\lambda}$ through $\alpha$ satisfies
$$
{\mathrm d} \alpha_{\lambda} +\frac{1}{2}[\alpha_{\lambda} \wedge
\alpha_{\lambda}]=0 \eqno(3.22)
$$
for every $\lambda$.
\end{Proposition}

It is clear that each $\alpha_{\lambda}$ is admissible whenever
$\alpha$ is, and $\alpha_{\lambda}$ generates the same loop. For
every $\alpha_{\lambda}$ satisfying (3.12), there exists a
one-parameter family of smooth maps $\Psi_{\lambda}:{\mathbb
D}\rightarrow G$ depending smoothly on $\lambda$ such that
$\Psi^{-1}_{\lambda} {\mathrm d}\Psi_{\lambda}= \alpha_{\lambda}$.
In  this paper, we shall normalize from here on $\Psi_{\lambda}$
by
$$
\Psi_{\lambda}(0,0) \equiv {\mathbf 1}. \eqno(3.23)
$$
Such normalized one-parameter family of maps $\Psi_{\lambda}$ is
called an {\it extended framing} of the harmonic map $\psi$. To
every harmonic maps $\psi:{\mathbb D}\rightarrow S^2_1$, there is
a naturally associated one-parameter family of harmonic maps
$\{\psi_{\lambda} \}$ such that $\psi_1=\psi$ parametrized by
$\lambda \in {\mathbb R}^{+}$. We shall therefore refer to
$\psi_{\lambda}$ as a {\it loop of harmonic maps} ({\it through\/}
$\psi$).

The extended framing $\Psi_\lambda$ of $\psi$ has its values in
${\tilde \Lambda}G_{\sigma}$ and can therefore be regarded as a
mapping $\Psi=\Psi_\lambda:{\mathbb D} \to {\tilde \Lambda}
G_{\sigma}$ into the twisted loop group. More generally we have
the following ({\it cf. \/} p.~116 in \cite{Gue}).

\begin{Proposition}{\rm (Harmonicity equation
in terms of extended framimgs)}
\newline
\noindent {\rm (1)} Let $\Psi=\Psi_\lambda: \mathbb{D}\to {\tilde
\Lambda}G_{\sigma}$ be a smooth map which satisfies the following
equations{\rm :}
$$
\Psi(0,0;\lambda)\equiv {\mathbf 1},
$$
$$
\Psi_\lambda^{-1}{\mathrm d}^{\prime} \Psi_\lambda= \mathrm{linear
\ in}\ \lambda \ (= A+\lambda B \ \mathrm{for}\ \mathrm{some}\  A,
B), \eqno(\Lambda_\sigma)
$$
$$
\Psi_\lambda^{-1} {\mathrm d}^{\prime \prime} \Psi_\lambda=
\mathrm{linear \ in}\ \lambda^{-1} \ (= C+\lambda^{-1} D \
\mathrm{for}\
 \mathrm{some}\  C, D)
$$
where $A$ and $C$ are $\mathfrak{k}$-valued and $B$ and $D$ are
$\mathfrak{m}$-valued $1$-forms.

Then $\psi_\lambda:=\mathrm{Ad}(\Psi_\lambda)\mathbf{k}^{\prime}$
defines a loop of harmonic maps into $S^2_1$.

\vspace{0.2cm}

\noindent {\rm (2)} Let $\psi:\mathbb{D}\to S^2_1$ be a harmonic
map. Then there exists a solution $\Psi_\lambda$ of
$(\Lambda_\sigma)$ such that
$\psi=\mathrm{Ad}(\Psi_1)\mathbf{k}^{\prime}$.

\end{Proposition}

\begin{Remark}
C.-H.~Gu \cite{Gu} investigated the Cauchy problem for Lorentzian
harmonic maps from ${\mathbb E}^2_1$ into $S^2_1$.
\end{Remark}

\section{ Weierstra{\ss}-type representation for
\newline harmonic maps }

{\large {\bf 4.1}}\hspace{0.15cm} In this section, we establish a
Weierstra{\ss}-type representation for \newline \noindent
Lorentzian harmonic maps into $S^2_1=G/K$.

\vspace{0.2cm}

We start by introducing the following linear spaces:

$$
\Lambda_{-\infty,1}:= \left \{ \xi^{\prime} \in {
\Lambda}{\mathfrak g}_\sigma \ | \
\xi^\prime=\sum_{k=-\infty}^{1}\xi^{\prime}_k\> \lambda^k \ \right
\}, \eqno(4.1)
$$
$$
\Lambda_{-1,\infty}:= \left \{ \xi^{\prime \prime} \in {
\Lambda}{\mathfrak g}_\sigma \ | \ \xi^{\prime \prime}=
\sum_{k=-1}^{\infty}\xi^{\prime \prime}_k\> \lambda^k \ \right \}.
\eqno(4.2)
$$

\begin{Definition}
Let us denote the space of all $\Lambda_{-\infty,1}$-valued
Lorentz holomorphic $1$-forms on $\mathbb{D}$ by
$\mathcal{P}^{\prime}$. Similarly we denote by
$\mathcal{P}^{\prime \prime}$ the space of
$\Lambda_{-1,\infty}$-valued Lorentz anti-holomorphic $1$-forms on
$\mathbb{D}$.
\end{Definition}

By definition the elements $\xi^{\prime} \in \mathcal{P}^\prime$
and $\xi^{\prime \prime} \in \mathcal{P}^{\prime \prime}$ have the
following form:
$$
\xi^\prime=\left( \sum_{k=-\infty}^{1}\xi^{\prime}_k(x)\>
\lambda^k \right ){\mathrm d}x,\ \ \ \xi^{\prime \prime}=\left(
\sum_{k=-1}^{\infty} \xi^{\prime \prime}_k(y)\> \lambda^k \right
){\mathrm d}y. \eqno(4.3)
$$
Here the coefficients $\xi^{\prime}_{k}(x)$ [resp. $\xi^{\prime
\prime}_{k}(y)$] are smooth functions of $x$ [resp. $y$]. We call
elements of $\mathcal{P}^{\prime}$ and $\mathcal{P}^{\prime
\prime}$ by $(1,0)$-{\it potentials} and $(0,1)$-{\it potentials}
for harmonic maps.

\vspace{0.2cm}

\noindent {\large {\bf 4.2}}\hspace{0.15cm} We will show below
that the space ${\mathcal P}^{\prime} \times {\mathcal P}^{\prime
\prime}$ can serve as the spaces of {\it Weierstra{\ss} data} for
the construction of harmonic maps from $\mathbb{D}$ into $S^2_1$.

\vspace{0.3cm}

\begin{Remark}
In the case of CMC surfaces in Euclidean $3$-space $\mathbb{E}^3$
the corresponding potentials are the {\it holomorphic} potentials
for harmonic Gau{\ss} maps from Riemann surfaces into $S^2$. In
particular for the description of finite type CMC surfaces in
${\mathbb E}^3$, it turns out to be very helpful to use these
potentials. See \cite{DH}.
\end{Remark}

\begin{Theorem}{\rm (Weierstra{\ss} type representation)}
\newline
\noindent Let $\{ \xi^{\prime}, \xi^{\prime \prime} \} \in
{\mathcal P}^{\prime} \times {\mathcal P}^{\prime \prime}$ be  a
potential.

Solve the two independent initial value problems {\rm :}
$$
{\mathrm d}^{\prime} \Psi^{\prime} =\Psi^{\prime}\ \xi^{\prime}, \
\ {\mathrm d}^{\prime \prime}\Psi^{\prime \prime} =\Psi^{\prime
\prime}\ \xi^{\prime \prime}, \eqno(4.4)
$$
$$
\Psi^{\prime}(x=0)=\Psi^{\prime \prime}(y=0)={\mathbf 1}.
\eqno(4.5)
$$
Then the Iwasawa decomposition
$$
(\Psi^{\prime},\Psi^{ \prime \prime})=
(\Psi,\Psi)(L_{-}^{-1},L_{+}^{-1}) \in \Delta ({\tilde
\Lambda}G_\sigma \times{\tilde \Lambda}G_\sigma) \cdot {\tilde
\Lambda}^{-}_{*}G_{\sigma} \times {\tilde \Lambda}^{+}G_{\sigma}
\eqno(4.6)
$$
gives an extended framing $\Psi$. Hence
$\psi:=\mathrm{Ad}(\Psi)\mathbf{k}^{\prime}$ gives a loop of
harmonic maps into $S^2_1$.
\end{Theorem}

\begin{pf}
\noindent From the Iwasawa splitting (Theorem 2.6) of
$(\Psi^{\prime},\Psi^{\prime \prime})$ we know locally around
$(0,0)$:
$$
\Delta ({\tilde \Lambda}G_\sigma \times {\tilde \Lambda}G_\sigma)
\ni (\Psi,\Psi)=(\Psi^{\prime},\Psi^{\prime \prime})
(L_{-},L_{+}).
$$
The Maurer-Cartan form $(\Psi,\Psi)^{-1} {\mathrm d}(\Psi,\Psi)$
is computed as
$$
(\Psi,\Psi)^{-1}\> {\mathrm d} (\Psi,\Psi) =(
L_{-}^{-1}\xi^{\prime}L_{-}+ L_{-}^{-1}{\mathrm d}L_{-},
L_{+}^{-1}\xi^{\prime \prime}L_{+}+ L_{+}^{-1}{\mathrm d}L_{+}),
\eqno(4.7)
$$
where we have used
$$
(\Psi^{\prime})^{-1} {\mathrm d} \Psi^{\prime}= \xi^{\prime}, \ \
(\Psi^{\prime \prime})^{-1} {\mathrm d} \Psi^{\prime \prime}=
\xi^{\prime \prime}.
$$

By the construction, the two components of $(4.7)$ are equal. We
call them $\alpha$. One reads off that the $\mathrm{d}x$-part of
the first component of the right hand side of $(4.7)$ only
involves positive powers of $\lambda$ with exponents equal to or
smaller than $1$, while the $\mathrm{d}x$-part of the second
component involves only exponents equal to or greater than $0$.
Thus in  $\alpha$ only the exponents $0$ and $1$ occur in the
$\mathrm{d}x$-part. A similar argument for the $\mathrm{d}y$-part
shows that it only involves the terms with exponents $-1$ and $0$.

Thus we obtain
$$
\alpha=\Psi^{-1} {\mathrm d}\Psi =\lambda^{-1}(\xi^{\prime
\prime}_{-1}{\mathrm d}y)+ (\xi^{\prime}_{0}{\mathrm d}x
+\xi^{\prime \prime}_{0}{\mathrm d}y)+ \lambda (\xi^{\prime}_{1}
{\mathrm d}x). \eqno(4.8)
$$
Hence $\Psi$ is an extended framing, since $\xi^{\prime}_{1},\>
\xi^{\prime \prime}_{-1} \in \mathfrak{m}, \ \ \xi^{\prime}_0,\>
\xi^{\prime \prime}_{0} \in \mathfrak{k}$. $\Box$
\end{pf}

\vspace{0.2cm}

\noindent {\large {\bf 4.3}}\hspace{0.15cm} This section is in
some sense the converse of the previous one. However, we shall
compute potentials only in some normal form, which, in general,
will involve singularities. A relation between smooth surfaces and
singularity free potentials of the form $(4.3)$ has not been
established yet.

\vspace{0.2cm}

Let $\psi:\mathbb{D}\to G/K$ be a harmonic map with an extended
framing $\Psi$. Then we perform both type of Birkhoff
decompositions (Theorem 2.3) for $\Psi$ as long as it is in both
big cells. First we consider the Birkhoff decomposition of the
type:
$$
{\tilde \Lambda}^{-}_{*}G_{\sigma} \times {\tilde
\Lambda}^{+}G_{\sigma} \subset {\tilde \Lambda}G_{\sigma}.
\eqno(4.9)
$$
Since the big cell $\mathscr{B}(-,+)$ of ${\tilde
\Lambda}G_\sigma$ is open and $\Psi$ is continuous, the set
$$
D_1:=\left \{ (x,y) \in \mathbb{D} \ \vert \ \Psi(x,y)\
\mathrm{belongs} \ \mathrm{to} \ \mathrm{the} \ \mathrm{big}\
\mathrm{cell}\ \mathscr{B}(-,+) \ \right \} \eqno(4.10)
$$
is an open subset of $\mathbb{D}$. Note that $D_1$ contains
$(0,0)$. Set $\mathcal{S}_1:=\mathbb{D}\setminus D_1$. Similarly,
we have $D_2$ and $\mathcal{S}_2$ for the splitting:
$$
{\tilde \Lambda}^{+}_{*}G_{\sigma} \times {\tilde
\Lambda}^{-}G_{\sigma} \subset {\tilde \Lambda}G_{\sigma}.
\eqno(4.11)
$$
We can perform the two Birkhoff splittings on the extended framing
$\Psi$ over $\mathbb{D}\setminus \mathcal{S},\ \ \mathcal{S}:=
\mathcal{S}_1 \cup \mathcal{S}_2$:
$$
\Psi=\Psi_{-}\>L_{+},\ \Psi_{-} \in {\tilde
\Lambda}_{*}^{-}G_{\sigma},\ L_{+}\in {\tilde
\Lambda}^{+}G_{\sigma}. \eqno(4.12)
$$
$$
\Psi=\Psi_{+}\>L_{-},\ \Psi_{+} \in {\tilde
\Lambda}_{*}^{+}G_{\sigma},\ L_{-}\in {\tilde
\Lambda}^{-}G_{\sigma}. \eqno(4.13)
$$
\noindent From the splitting (4.12), we obtain
$$
\Psi^{-1}_{-}{\mathrm d}\Psi_{-}= L_{+}(\Psi^{-1}{\mathrm
d}\Psi)L_{+}^{-1}- \mathrm{d}L_{+}\>L_{+}^{-1}. \eqno(4.14)
$$
Since $\alpha_\lambda =\alpha_{0} +\lambda \alpha^{\prime}_{1}
+\lambda^{-1}\alpha^{\prime \prime}_{1}$, we have
$$
L_{+}(\Psi^{-1}{\mathrm d}\Psi)L_{+}^{-1} =L_{+} (\alpha_{0}
+\lambda \alpha^{\prime}_{1} +\lambda^{-1}\alpha^{\prime
\prime}_{1}) L_{+}^{-1}. \eqno(4.15)
$$
Since $L_{+} \in {\tilde \Lambda}^{+}G_{\sigma}$, $L_{+}$ has the
decomposition
$$
L_{+}=\sum_{k\geq 0}L^{+}_{k}(x,y)\> \lambda^k. \eqno(4.16)
$$
Compare the left and right hand sides of the $(1,0)$-part of
$(4.14)$.
\newline
The left hand side of $(4.14)$ contains negative powers of
$\lambda$. On the other hand, the right hand side contains
nonnegative powers of $\lambda$ only. Thus we have
$$
\Psi_{-}^{-1}\frac{\partial \Psi_{-}} {\partial x}=0. \eqno(4.17)
$$
Therefore $\Psi_{-}$ depends only on $y$ and $\lambda$.

Next by comparing the left and right hand sides of the
$(0,1)$-part of $(4.14)$, we obtain:
$$
\Psi_{-}^{-1}\> {\mathrm d}^{\prime \prime} \Psi_{-}= \lambda^{-1}
\left \{ \> \mathrm{Ad} \left(\> L^{+}_{0} \> \right ) \>
\alpha_1^{\prime \prime} \right \}{\mathrm d}y. \eqno(4.18)
$$
The left hand side of $(4.14)$ depends only on $y$ (and
$\lambda$). Hence
$$
\eta^{\prime \prime}:= \mathrm{Ad} \left(\> L^{+}_{0} \> \right )
\> \alpha_1^{\prime \prime} \> {\mathrm d}y \eqno(4.19)
$$
is an ${\mathfrak m}$-valued anti-holomorphic $1$-form on
$\mathbb{D}\setminus \mathcal{S}$.

By similar arguments for $\Psi_{+}=\Psi \ L_{-}^{-1}$, where
$L_{-}=\sum_{k\leq 0}L^{-}_{k}\lambda^k$, we have
$$
\Psi_{+}^{-1}\frac{\partial \Psi_{+}} {\partial y}=0. \eqno(4.20)
$$
Thus $\Psi_{+}$ depends only on $x$ (and $\lambda$). Moreover the
$1$-form
$$
\eta^{\prime}:= \mathrm{Ad} \left(\> L^{-}_{0} \> \right ) \>
\alpha_1^{\prime} \> {\mathrm d}x \eqno(4.21)
$$
is an ${\mathfrak m}$-valued  holomorphic $1$-form on
$\mathbb{D}\setminus \mathcal{S}$. Obviously
$$
\xi^{\prime}:=\lambda \eta^{\prime} \in {\mathcal P}^{\prime},\ \
\xi^{\prime \prime}:=\lambda^{-1} \eta^{\prime \prime} \in
{\mathcal P}^{\prime \prime}. \eqno(4.22)
$$
One can check that the pair of potentials $\xi^{\prime},\
\xi^{\prime \prime}$ reproduces the harmonic map $\psi$ via
Weierstra{\ss} representation. ({\it cf.\/} Lemma 4.5 and Theorem
4.10 in \cite{DPW} and Theorem 2.1 in \cite{Wu}.)

\begin{Theorem}{\rm(The normalized potentials)}
\newline
\noindent Let  $\psi:\mathbb{D}\to G/K$ be a harmonic map with
$\psi(0,0)=K$ and $\Psi$ an extended framing of $\psi$. Then there
exists an open subset $(0,0) \in \mathbb{D}\setminus{\mathcal S}$
on which $\Psi$ splits into
$$
\Psi=\Psi_{-}L_{+}=\Psi_{+}L{-}, \eqno(4.23)
$$
$$
\Psi_{\pm}\in {\tilde \Lambda}^{\pm}_{*}G_{\sigma},\ \ L_{\pm} \in
{\tilde \Lambda}^{\pm}G_{\sigma}.
$$
The $1$-forms $\eta^\prime$, $\eta^{\prime \prime}$ defined by
$$
\eta^{\prime}(x)=\Psi_{+}^{-1} {\mathrm d}\> \Psi_{+}\> \lambda,\
\ \ \eta^{\prime \prime}(y)=\Psi_{-}^{-1} {\mathrm d}\> \Psi_{-}\>
\lambda^{-1} \eqno(4.24)
$$
are an $\mathfrak{m}$-valued holomorphic $1$-form and an
anti-holomorphic $1$-form on $\mathbb{D}\setminus{\mathcal S}$
respectively.

\vspace{0.2cm}

Conversely, any harmonic map $\psi:\mathbb{D}\to G/K$ with
$\psi(0,0)=K$ can be constructed from a pair of
$\mathfrak{m}$-valued $1$-forms $\eta^\prime,\ \eta^{\prime
\prime}$ which are holomorphic and anti-holomorphic respectively.
The harmonic map $\psi$ is constructed via the Weierstra{\ss}
representation with potentials
$$
\xi^{\prime}:=\lambda \eta^{\prime},\ \ \xi^{\prime
\prime}:={\lambda}^{-1} \eta^{\prime \prime}. \eqno(4.25)
$$
\end{Theorem}

\vspace{0.2cm}

The pair of $1$-forms $\{\eta^\prime,\ \eta^{\prime \prime}\}$ is
defined uniquely.
\newline
Following Wu \cite{Wu} and \cite{T}, we call the pair
$\{\eta^\prime,\ \eta^{\prime \prime}\}$ (or the pair
$\{\xi^{\prime},\> \xi^{\prime \prime} \}$ defined by
$\xi^{\prime}= \lambda \eta^{\prime},\ \xi^{\prime
\prime}:=\lambda^{-1} \eta^{\prime \prime}$ ) the {\it normalized
potentials} for $\psi$ with the origin as the reference point.

\smallskip

 Up to now
 we have only very little information
 on the singular set.
 Thus in this paper,
 we restrict
 our attention to holomorphic and
 anti-holomorphic
 potentials
 $\{\xi^{\prime}, \xi^{\prime \prime}\}$.

\begin{Remark}
In Euclidean CMC surface geometry, the holomorphic potentials can
be extended meromorphically to the simply connected Riemann
surface $\mathbb{D}$. The poles of the potentials are in the
singular set $\mathcal{S}$.
\end{Remark}

\begin{Remark}
Balan and the first named author studied Weierstra{\ss}-type
representation of harmonic maps from Riemann surfaces into
noncompact Riemannian symmetric spaces \cite{BD2}.
\end{Remark}


\section{Normalized potentials
for timelike \newline CMC surfaces}

\noindent {\large {\bf 5.1}}\hspace{0.15cm} In this section, we
shall apply the Weierstra{\ss} representation for harmonic maps
into $S^2_1$ to timelike CMC surfaces.

\vspace{0.2cm}

Let $\varphi:{\mathbb D}\rightarrow {\mathbb E}^3_1$ be a timelike
CMC surface with Gauss map $\psi$ and coordinate frame ${\hat
\Phi}$ (See (1.3)). The Gauss-Codazzi equations (G) and (C) of
$\varphi$ are invariant under the deformation:
$$
Q \ \longmapsto Q_{\lambda}:=\lambda Q,\ \ R \ \longmapsto
R_{\lambda}:= \lambda^{-1} R, \ \lambda \in {\mathbb R}^{+}.
\eqno(5.1)
$$
Integrating the Frenet equations $(1.4)$ with $Q_{\lambda}$ and
$R_{\lambda}$, one obtains a one-parameter family of timelike
surfaces $\{ {\hat \varphi} _{\lambda} \}$. This deformation does
not effect the induced metric and the mean curvature. Hence all
the surfaces $\{ {\hat \varphi}_{\lambda} \}$ are isometric and
have the same constant mean curvature. The family $\{{\hat
\varphi}_\lambda\}$ is called the {\it associated family} of
$\varphi$.

This one-parameter deformation of ${\hat \varphi}$ satisfies the
following Lax equations:
$$
\frac{\partial}{\partial x} {\hat \Phi}_{\lambda}={\hat
\Phi}_{\lambda} {\hat U}(\lambda),\ \ \frac{\partial}{\partial y}
{\hat \Phi}_{\lambda}={\hat \Phi}_{\lambda} {\hat V}(\lambda),\ \
\eqno(5.2)
$$
$$
{\hat U}(\lambda) =\left (
\begin{array}{cc}
-\frac{1}{4}\omega_{x} & -\lambda
Qe^{-\frac{\omega}{2}} \\
\frac{H}{2}e^{\frac{\omega}{2}} & \frac{1}{4}\omega_x
\end{array}\right )
,\ \ {\hat V}(\lambda) =\left (
\begin{array}{cc}
\frac{1}{4}\omega_y &
-\frac{H}{2}e^{\frac{\omega}{2}} \\
\lambda^{-1}Re^{-\frac{\omega}{2}} & -\frac{1}{4}\omega_y
\end{array}
\right ). \eqno(5.3)
$$

A solution ${\hat \Phi}_\lambda$ to $(5.2)$ is not an extended
framing for the harmonic Gau{\ss} map $\psi$ in the sense of the
definition given in Section 3.5, since the $\lambda$-distribution
does not fit.

To relate these two $S^1$-families to each other we perform the
transformation
$$
\Phi:= g(\lambda)^{-1}\>{ {\hat \Phi}_{\lambda^2}}\> g(\lambda), \
\ \lambda \in \mathbb{R^+} \eqno(5.4)
$$
with
$$
g(\lambda)= \left (
\begin{array}{cc}
 0 & -\sqrt{\lambda} \\
1/\sqrt{\lambda} & 0
\end{array}
\right ) .
$$

Then the coefficient matrices of the Lax pair $\{{\hat U},{\hat
V}\}$ are changing into
$$
U=g(\lambda)^{-1}{{\hat U}(\lambda^2)}g(\lambda)= \left (
\begin{array}{cc}
\frac{1}{4}\omega_{u} &
-\lambda\frac{H}{2}e^{\frac{\omega}{2}} \\
\lambda Qe^{-\frac{\omega}{2}}  & -\frac{1}{4}\omega_u
\end{array}
\right ), \eqno(5.5)
$$
$$
V=g(\lambda)^{-1} {{\hat V}(\lambda^2)}g(\lambda)= \left (
\begin{array}{cc}
-\frac{1}{4}\omega_{v} &
-\lambda^{-1}Re^{-\frac{\omega}{2}} \\
\lambda^{-1} \frac{H}{2}e^{-\frac{\omega}{2}}  &
\frac{1}{4}\omega_v
\end{array}
\right ).
$$

Comparing with \cite{I} we conclude

\begin{Proposition}
{\rm (}Sym formula{\rm )}

Let $\Phi$ be a solution to the Lax equations {\rm (5.4) and
(5.5)}. Then
$$
\varphi_{\lambda}=-\frac{1}{H} \left \{ \frac{\partial}{\partial
t} \Phi \cdot \Phi^{-1} +\frac{1}{2} \mathrm{Ad} (\Phi){\mathbf k}
^{\prime} \right \}, \ \ \lambda=e^t \in \mathbb{R}^+
\eqno({\mathrm S})
$$
describes a real loop of timelike surfaces of constant mean
curvature $H$. The first fundamental form $I$ and the Gau{\ss} map
$N_\lambda$ of $\varphi_\lambda$ are given by
$$
I=e^\omega {\mathrm d}x{\mathrm d}y \ \ and \ \
N_{\lambda}=\mathrm{Ad}(\Phi) {\mathbf k}^{\prime}. \eqno(5.6)
$$
The logarithmic derivative part $\varphi^{K}_{\lambda}$ of
$\varphi_\lambda${\rm :}
$$
\varphi^{K}_{\lambda}=-\frac{1}{H} \frac{\partial}{\partial t}\Phi
\cdot \Phi^{-1}
$$
describes a real loop of timelike surfaces with constant
Gau{\ss}ian curvature $4H^2$.

\end{Proposition}

\vspace{0.2cm}

\noindent {\large {\bf 5.2}}\hspace{0.15cm} Next we calculate the
normalized potential $\{\xi^{\prime},\> \xi^{\prime \prime}\}$ in
terms of the fundamental quantites of the timelike CMC surface
$\varphi$. More precisely we shall clarify the role of the
normalized potentials in the construction of timelike CMC
surfaces.

\vspace{0.2cm}

First decomopose $\alpha:=\Phi^{-1}{\mathrm d}\Phi$ as
$$
\Phi^{-1}{\mathrm d}\Phi=\alpha_{0} +\lambda
\>\alpha_{1}^{\prime}+ \lambda^{-1}\>\alpha_{1}^{\prime \prime}.
$$
Next we perform  the Birkhoff decompositions of $\Phi$ with regard
to
$$
{\tilde \Lambda}^{+}_{*}G_\sigma \times {\tilde
\Lambda}^{-}G_\sigma \subset {\tilde \Lambda}G_{\sigma},\ \
{\tilde \Lambda}^{-}_{*}G_\sigma \times {\tilde
\Lambda}^{+}G_\sigma \subset {\tilde \Lambda}G_\sigma,\ \
$$
over the big cells $\mathscr{B}(-,+)$ and $\mathscr{B}(+,-)$ of
${\tilde \Lambda}G_\sigma$:
$$
\Phi=\Phi_{+}\>L_{-} =\Phi_{-}\>L_{+}. \eqno(5.7)
$$

We have seen above that with
$$
L_{+}=\sum_{k\geq 0}L^{+}_{k}\lambda^k,\ \ L_{-}=\sum_{k\leq
0}L^{-}_{k}\lambda^k, \eqno(5.8)
$$
the potentials
$$
\xi^{\prime}:=\Phi_{+}^{-1}\> {\mathrm d}^{\prime} \Phi_{+},\ \
\xi^{\prime \prime}: =\Phi_{-}^{-1}\> {\mathrm d}^{\prime \prime}
\Phi_{-} \eqno(5.9)
$$
have the  form:

$$
\xi^{\prime}= \lambda \ \left \{
L^{+}_{0}(x,y)\>\alpha^{\prime}_{1} \> L^{+}_{0}(x,y)^{-1} \right
\}, \eqno(5.10)
$$
$$
\xi^{\prime \prime}= \lambda^{-1} \ \left \{
L^{-}_{0}(x,y)\>\alpha^{\prime \prime}_{1} \> L^{-}_{0}(x,y)^{-1}
\right \}. \eqno(5.11)
$$

More explicitly we have
$$
\Phi_{+}^{-1} \frac{\partial \Phi_{+}} {\partial x}= \lambda \left
\{ L^{-}_{0} \left(
\begin{array}{cc}
0 & -\frac{H}{2}e^{\omega/2} \\
Q(x)e^{-\omega/2}  & 0
\end{array}
\right ) \ (L^{-}_{0})^{-1} \right \}, \eqno(5.12)
$$
$$
\Phi_{-}^{-1} \frac{\partial \Phi_{-}} {\partial y}= \lambda^{-1}
\left \{ L^{+}_{0} \left(
\begin{array}{cc}
0 & -R(y)e^{-\omega/2} \\
\frac{H}{2}e^{\omega/2} & 0
\end{array}
\right ) \ (L^{+}_{0})^{-1} \right \}. \eqno(5.13)
$$
Comparing  the two sides of $(5.13)$ we see that the
 left hand side depends on $y$ and $\lambda$ only, while
the right hand side depends on $x$, $y$ and $\lambda$. Recall from
the Coddazi equation (C) that since the mean curvature  $H$ is
constant, $R$ depends only on $y$. Hence we derive
$$
\Phi^{-1}_{-} \frac{\partial \Phi_{-}} {\partial y}= \lambda^{-1}
\left \{ L^{+}_{0}(0,y) \left (
\begin{array}{cc}
0 & -R(y)e^{-\omega(0,y)/2} \\
\frac{H}{2}e^{\omega(0,y)/2} & 0
\end{array}
\right ) L^{+}_{0}(0,y)^{-1} \right \}.
$$

We abbreviate $L^{+}_{0}(0,y)$ by $W_{0}(y)$ and  compute
$W_{0}(y)$ more explicitly. \noindent From (5.7) we obtain the
Maurer-Cartan equation:
$$
L_{+} \> {\alpha} \> L_{+}^{-1} - {\mathrm d}L_{+}\> L_{+}^{-1} =
\Phi_{-}^{-1} {\mathrm d}{\Phi}_{-}. \eqno(5.14)
$$

Comparing the $\lambda^0$-terms of both sides in $(5.14)$ we
obtain,
$$
0=W_{0}\ \alpha_{0}^{\prime \prime}\ W_{0}^{-1}- {\mathrm
d}W_{0}\>W_{0}^{-1}. \eqno(5.15)
$$
Hence
$$
{\mathrm d}W_{0}=W_{0}\> \alpha_{0}^{\prime \prime}. \eqno(5.16)
$$
Namely we have
$$
\frac{{\mathrm d}W_{0}}{\mathrm{d}y}= W_{0} \left \{ \frac
{\omega_{y}(0,y)}{4} \right \}{\mathbf k}^{\prime}. \eqno(5.17)
$$
Thus $W_{0}(y)$ is given explicitly by
$$
W_{0}(y)= \left (
\begin{array}{cc}
e^{-\omega(0,y)/4+c_1} & 0 \\
0 & e^{\omega(0,y)/4+c_2}
\end{array}
\right ),\ \ c_1,\ c_2 \in \mathbf{R}. \eqno(5.18)
$$
Since we require the initial condition:
$$
W_{0}(0)=\mathbf{1},
$$
we see that  the $(0,1)$-potential $\xi^{\prime \prime}$ is given
by

$$
\xi^{\prime \prime}=\lambda^{-1} \left (
\begin{array}{cc}
0 & -R(y) e^
{\{-\omega(0,y)+\omega (0,0)/2\}} \\
\frac{H}{2}e^ {\{\omega(0,y)-\omega (0,0)/2\}} & 0
\end{array}
\right ) {\mathrm d}y. \eqno(5.19)
$$

By similar arguments for $\Phi_{+}=\Phi \ L_{-}^{-1}$, we obtain
$$
\Phi^{-1}_{+} \frac{\partial \Phi_{+}} {\partial x}= \lambda \left
\{ L^{-}_{0}(x,0) \left (
\begin{array}{cc}
0 &
-\frac{H}{2}e^{\omega(x,0)/2} \\
Q(x)e^{-\omega(x,0)/2}& 0 \
\end{array}
\right ) L^{-}_{0}(x,0)^{-1} \right \}.
$$
The function $\Gamma_{0}(x):=L^{-}_{0}(x,0)$ is a solution to
$$
\frac{\mathrm{d}\Gamma_{0}}{\mathrm{d}x}= \Gamma_{0} \left \{
\frac {-\omega_{x}(x,0)}{4} \right \}{\mathbf k}^{\prime}.
\eqno(5.20)
$$
Thus we have
$$
\Gamma_{0}(x)= \left (
\begin{array}{cc}
e^{\omega(x,0)/4+c_3} & 0 \\
0 & e^{-\omega(x,0)/4+c_4}
\end{array}
\right ),\ \ c_3,\ c_4 \in \mathbb{R}. \eqno(5.21)
$$
Since we require the initial condition
$$
\Gamma_{0}(0,0)=\mathbf{1},
$$
the $(1,0)$-potential $\xi^{\prime}$ is given by
$$
\xi^{\prime}=\lambda \left (
\begin{array}{cc}
0 & -\frac{H}{2}e^
{\{\omega(x,0)-\omega (0,0)/2\}} \\
Q(x) e^ {\{-\omega(x,0)+\omega (0,0)/2\}} & 0
\end{array}
\right ) {\mathrm d}x. \eqno(5.22)
$$

\begin{Theorem}
The normalized potentials with the origin as a reference point are
given by
\[
\xi^{\prime}=\lambda \left (
\begin{array}{cc}
0 & -\frac{H}{2}
f(x) \\
Q(x) /f(x) & 0
\end{array}
\right ) {\mathrm d}x,
\]
\[
\xi^{\prime \prime}=\lambda^{-1} \left (
\begin{array}{cc}
0 &
-R(y)/g(y) \\
\frac{H}{2}g(y) & 0
\end{array}
\right ) {\mathrm d}y.
\]
Here the functions $f(x)$ and $g(y)$ are given by
$$
f(x)=\mathrm{exp} \> \{\omega(x,0)-\omega (0,0)/2\}, \ \
g(y)=\mathrm{exp} \> \{\omega(0,y)-\omega (0,0)/2\}.
$$
\end{Theorem}

\begin{Remark} From this result one sees that the
normalized potentials are essentially the Hopf differentials
together with the values of the conformal factor of the metric
restricted to one pair of null coordinate lines.
\end{Remark}

\begin{Remark}{\rm (Nonlinear
d'Alembert's formulas)}

It is known that every solution for the linear wave equation:
$$
\omega_{xy}=0
$$
can be written as the sum of a Lorentz holomorphic function and a
Lorentz anti-holomorphic function by the so called {\it d'Alembert
formula}:
$$
\omega(x,y)=f(x)+g(y).
$$
The Weierstra{\ss}-type representation (Theorem 4.3) toghether
with Theorem 5.2 provide us with nonlinear analogues of
d'Alembert's formula for the hyperbolic sinh-Gordon equation, the
Liouville equation and the hyperbolic cosh-Gordon equation. More
precisely for any initial data $f(x)$ and $g(y)$, the normalized
potentials with $Q=H/2,\ R=\epsilon H/2,\ \epsilon =0,\pm 1$
produce solutions $\omega$ to the hyperbolic sinh-Gordon equation
($\epsilon=1$), the Liouville equation ($\epsilon=0$) and the
hyperbolic cosh-Gordon equation ($\epsilon=-1$) via the
Weierstra{\ss}-type representation (Theorem 4.3). (See Proposition
1.3, Proposition 1.5 and Remark 1.6.)
\end{Remark}


\section{Examples}

Even though the results of the previous sections give a one-to-one
correspondence between special potentials and timelike CMC
surfaces, making this correspondence explicit is a different
matter.
\par
To give examples we start with the easiest case from the potential
point of view.

\vspace{0.2cm} \noindent
{\large {\bf 6.1}} (Cylinders and pseudospheres)
\begin{Example}{\rm (}hyperbolic cylinders{\rm )}
Let us take the following potentials:
$$
\xi^{\prime}=\frac{\lambda}{4} \left (
\begin{array}{cc}
0 & -1 \\
-1 & 0
\end{array}
\right ) {\mathrm d}x,\ \ \ \xi^{\prime \prime}=
\frac{\lambda^{-1}}{4} \left (
\begin{array}{cc}
0 & 1 \\
1 & 0
\end{array}
\right) {\mathrm d}y. \eqno(6.1)
$$
Then the normalized potentials $\{\xi^{\prime},\ \xi^{\prime
\prime}\}$ produce the associated family of timelike CMC surface
$\varphi:{\mathbb R}^2 \to \mathbb{E}^3_1$ of mean curvature
$1/2$;
$$
\varphi(x,y)=(\sinh \frac{x-y}{2}, \frac{x+y}{2}, \cosh
\frac{x-y}{2}). \eqno(6.2)
$$
The image of $\varphi$ is a timelike {\it hyperbolic} cylinder.
The induced metric of ${\mathbb R}^2$ is $I={\mathrm d}x{\mathrm
d}y$. Thus $\varphi$ is an isometric imbedding of the Minkowski
plane $({\mathbb R}^2,{\mathrm d}x{\mathrm d}y)$ into
$\mathbb{E}^3_1$.
\end{Example}

\begin{Example}
{\rm (}Circular cylinders{\rm )}

The normalized potentials $\{\xi^{\prime},\> \xi^{\prime
\prime}\}$ defined by

$$
\xi^{\prime}=\frac{\lambda}{4} \left (
\begin{array}{cc}
0 & -1 \\
1 & 0
\end{array}
\right ) {\mathrm d}x,\ \ \ \xi^{\prime
\prime}=\frac{\lambda^{-1}}{4} \left (
\begin{array}{cc}
0 & -1 \\
1 & 0
\end{array}
\right ) {\mathrm d}y \eqno(6.3)
$$
produce the associated family of timelike CMC surface
$\varphi:{\mathbb R}^2 \to \mathbb{E}^3_1$ of mean curvature
$1/2$;
$$
\varphi(x,y)=(\frac{x-y}{2},\sin \frac{x+y}{2}, -\cos
\frac{x+y}{2}). \eqno(6.4)
$$
The image of $\varphi$ is a timelike {\it circular} cylinder. The
induced metric of ${\mathbb R}^2$ is $I={\mathrm d}x{\mathrm d}y$.
Thus $\varphi$ is an isometric immersion of the Minkowski plane
$({\mathbb R}^2,{\mathrm d}x{\mathrm d}y)$ into $\mathbb{E}^3_1$.
\end{Example}
Note that both, timelike hyperbolic cylinder and timelike circular
cylinder, correspond to the {\it vacuum solution} of the
hyperbolic sinh-Gordon equation.

\begin{Example}(Totally umbilical pseudosphere)

The normalized potentials $\{\xi^{\prime},\>\xi^{\prime \prime}\}$
with $Q=R=0$ produce totally umbilical pseudosphere in
$\mathbb{E}^3_1$.

\end{Example}

\vspace{0.2cm}

\noindent {\large {\bf 6.2}}\hspace{0.15cm} Next we shall give
examples of timelike CMC surfaces with repeated real principal
curvatures. Namely, timelike CMC surfaces with $QR=0$. Such
surfaces have no Euclidean counterparts.

As in the case of indefinite affine spheres with $AB=0$ (See
Section 9.2 in  \cite{DE}), timelike CMC surfaces with $QR=0$ have
specific shapes. To investigate such surfaces, we recall the
notion of a $B$-scroll introduced by L.~Graves \cite{G}.
\begin{Definition}
Let $\gamma=\gamma(s)$ be a smooth curve in ${\mathbb E}^3_1$
defined on an interval $I\subset \mathbb{R}$. Then $\gamma$ is
said to be a {\it null Frenet curve} if \vspace{0.2cm}
\newline
(1) $\langle \gamma^{\prime}, \gamma^{\prime} \rangle=0$,
\vspace{0.2cm}
\newline
(2) there exist vector fields $A,B,C$ along $\gamma$ and two
functions $\kappa$ and $\tau$ such that
$$
A=\gamma^{\prime},\ \langle A,B\rangle=1,
$$
$$
\langle A,A \rangle= \langle B,B \rangle=0,\ \ \langle C,C
\rangle=1,
$$
$$
\langle A,C \rangle= \langle B,C \rangle=0. \eqno(6.5)
$$
\[
\frac{\mathrm d}{{\mathrm d}s}(A,B,C)= (A,B,C) \left (
\begin{array}{ccc}
0 & 0 & -\tau \\
0 & 0 & -\kappa \\
\kappa & \tau & 0
\end{array}
\right ).
\]
The frame fieled $L=(A,B,C)$ along $\gamma$ is called the {\it
null Frenet frame field} of $\gamma$. The two functions $\kappa$
and $\tau$ are called the {\it curvature} and {\it torsion} of
$\gamma$ respectively.
\end{Definition}

\begin{Definition}
Let $\gamma$ be a null Frenet curve in $\mathbb{E}^3_1$. The ruled
surface
$$
\varphi(s,t)=\gamma(s)+tB(s): I\times \mathbb{R}^{*}\rightarrow
 \mathbb{E}^3_1
\eqno(6.6)
$$
is called $B$-{\it scroll} of $\gamma$.
\end{Definition}
Since
$$
I=(t\tau)^2{\mathrm d}s^2+ 2{\mathrm d}s\,{\mathrm d}t,\ \
$$
every $B$-scroll is timelike. (In fact $\det I=-1$.) The mean
curvature of $\varphi$ is $H(s,t)=\tau(s)$. Thus $\varphi$ is of
constant mean curvature if and only if $\tau$ is constant.

\medskip

A specific example of a null Frenet curve with constant torsion
$1$ is
$$
\gamma(s)=\left( \frac{\sinh (2s)}{2}, \frac{\cosh (2s)}{2}, s
\right)
$$
with null Frenet frame field:
$$
A(s)=\left( \cosh 2s, \sinh 2s, 1 \right),
$$
$$
B(s)=\frac{1}{2} \left( -\cosh (2s), -\sinh (2s),1 \right),
$$
$$
C(s)= \left(-\cosh 2s, -\sinh 2s, 0 \right ).
$$
The $B$-scroll $\varphi$ of $\gamma$:
$$
\varphi(s,t)= \left( \frac{\sinh (2s)-t\cosh(2s)}{2}, \frac{\cosh
(2s)-t\sinh(2s)}{2}, s+\frac{t}{2} \right )
$$
is a timelike surface with constant mean curvature $1$.
(\cite{Mc}, p.~33.)

\medskip
Now let $\varphi$ be a constant mean curvature $B$-scroll of a
null Frenet curve $\gamma$. The local coordinate system $(x,y)$
defined by
$$
x:=s-\frac{2}{\tau^{2}t}-\frac{1}{2},\ \ y:=s+\frac{1}{2}.
\eqno(6.7)
$$
is a null coordinate system of $\varphi$. Then
$$
I=e^{\omega(x,y)}\> {\mathrm d}x{\mathrm d}y,\ \
e^{\omega(x,y)}=\frac{4}{H^2\{(x-y)+1\}^2}, \eqno(6.8)
$$
$$
Q(x)=0,\ R(y)=H(y-\frac{1}{2}). \eqno(6.9)
$$
Thus the normalized potentials are given by

$$
\xi^{\prime}=\lambda \left (
\begin{array}{cc}
0 & (x+1)^{-2} \\
0 & 0
\end{array}
\right){\mathrm d}x,
$$
$$
\xi^{\prime \prime}= \lambda^{-1} \left (
\begin{array}{cc}
0 & -R(y)H{(1-y)^2}/2 \\
(1-y)^{-2} & 0
\end{array}
\right){\mathrm d}y.
$$

\vspace{0.5cm}

Conversely, we shall prove that timelike CMC surfaces derived from
normalized potentials of the form:
$$
\xi^{\prime}=\lambda \left (
\begin{array}{cc}
0 & -\frac{H}{2}f(x) \\
0 & 0
\end{array}
\right ) {\mathrm d}x,\ \ \ \xi^{\prime \prime}= \lambda^{-1}
\left(
\begin{array}{cc}
0 &
-R(y)/g(y) \\
\frac{H}{2}g(y) & 0
\end{array}
\right ) {\mathrm d}y. \eqno(6.10)
$$
are $B$-scrolls of constant mean curvature $H$.

\vspace{0.2cm}

To this end we consider the initial value problems:
$$
{\mathrm d}^{\prime}\Phi_{+}= \Phi_{+}\>\xi^{\prime},\ \ {\mathrm
d}^{\prime \prime} \Phi_{-}= \Phi_{-}\>\xi^{\prime \prime},
$$
$$
\Phi_{+}(x=0)= \Phi_{-}(y=0)=\mathbf{1}
$$
with potential $(6.10)$.
\par
The solution $\Phi_{+}$ is easily obtained as
$$
\Phi_{+}= \left (
\begin{array}{cc}
1 & -\lambda F(x) \\
0 & 1
\end{array}
\right),\ \ F(x)=\frac{H}{2}\int^x_0 f(x){\mathrm d}x. \eqno(6.11)
$$
The inverse loop of $\Phi_{+}$ is
$$
\Phi_{+}^{-1}= \left(
\begin{array}{cc}
 1 & \lambda F(x) \\
0 & 1
\end{array}
\right). \eqno(6.12)
$$
We shall get the extended framing $\Phi$ derived from
$\{\xi^{\prime},\ \xi^{\prime \prime}\}$ by the technique used in
\cite{DH}. (See section 3.6 of \cite{DH}). Express $\Phi_{-}$ as
$$
\Phi_{-}= \left(
\begin{array}{cc}
a & b \\ c& d
\end{array}
\right )= \mathbf{1}+ \sum_{j=1}^{\infty} \left(
\begin{array}{cc}
a_i & b_i \\
c_{i} & d_{i}
\end{array}
\right )\lambda^{-j}\\. \eqno(6.13)
$$
Take a ${\tilde \Lambda}G^{+}_{\sigma}$-valued map $\mathscr{G}$;
$$
\mathscr{G}(x,y)= \left(
\begin{array}{cc}
1 & -\lambda F/(1+c_1 F)  \\
0 & 1
\end{array}
\right ). \eqno(6.14)
$$
Then $ {\tilde \Phi}: =\Phi_{+}^{-1}\> \Phi_{-} \> \mathscr{G} $
is computed as
$$
{\tilde \Phi} ={\tilde \Phi}_{0}+ \sum_{j\geq 1}{\tilde
\Phi}_{j}\lambda^{-j}= \left (
\begin{array}{cc}
1+c_1 F & 0 \\
0 & 1/(1+c_1 F)
\end{array}
\right ) +\mathrm{negative} \ \mathrm{powers}.
\]

Put ${\tilde \Phi}_{+} ={\tilde \Phi}_{0}\>\mathscr{G}^{-1}$ and
${\tilde \Phi}_{-} :={\tilde \Phi}\ {\tilde \Phi}_{0}^{-1}$. Then
${\tilde \Phi}_{-} \in {\tilde \Lambda}^{-}G_{\sigma},\ {\tilde
\Phi}_{+} \in {\tilde \Lambda}_{*}^{+}G_{\sigma}$ and $
{\Phi_+}^{-1} \Phi_{-} = {\tilde \Phi}_{-} {\tilde \Phi}_{+}$.
Thus we obtain

$$
{\tilde \Phi}_{-}=L_{-} \eqno(6.15)
$$
and  the extended framing $\Phi$ is computed as follows:
$$
\Phi=\Phi_{+}L_{-}= \Phi_{+}{\tilde \Phi}_{-}= \Phi_{-}{\tilde
\Phi}_{+}^{-1}. \eqno(6.16)
$$

By the Sym formula, Proposition 5.1, we obtain the constant mean
curvature immersion
$$
\varphi(x,y)= \gamma(y)+q(x,y)B(y), \eqno(6.18)
$$
$$
\gamma(y):=-\frac{1}{H} \left \{ \frac{\partial \Phi_-} {\partial
t}\Phi_-^{-1} +\frac{1}{2} \mathrm{Ad}(\Phi_-)\mathbf{k}^\prime
\right \}, \eqno(6.19)
$$
$$
B(y):= \frac{\lambda}{H^2} \mathrm{Ad}(\Phi_-)
(\mathbf{j}^{\prime}-\mathbf{i}), \ \
q(x,y)=\frac{F(x)}{H\{1+c_1(y)F(x)\}}. \eqno(6.20)
$$
A direct computation shows that $\gamma(y)$ is a null Frenet curve
and $\varphi$ is the $B$-scroll of $\gamma$.

\bigskip
({\sc Dorfmeister}) \quad {\sc Fakult{\"a}t f{\"u}r Mathematik,
Technische Universit{\"a}t M{\"u}nchen, Arcisstrasse 21, D-80290,
M{\"u}nchen, Germany }

\smallskip

{\it E-mail adress}: {\tt dorfm@ma.tum.de}

\medskip
({\sc Inoguchi}) \quad {\sc Department of Applied Mathematics,
\newline
Fukuoka University, Fukuoka, 814-0180, Japan}

\smallskip

{\it E-mail adress}: {\tt inoguchi@bach.sm.fukuoka-u.ac.jp}

\medskip
({\sc Toda}) \quad {\sc Department of Mathematics and Statistics,
\newline Texas Tech University, Lubbock, TX, 79409-1042, U.S.A.}

\smallskip

{\it E-mail adress}: {\tt mtoda@math.ttu.edu}

\end{document}